\documentclass[aos,MSNbibl,seceqn,nameyear,dvips]{arximspdf}
\usepackage{graphicx}
\doi{10.1214/11-AOS962}
\volume{40}
\issue{1}
\pubyear{2012}
\firstpage{353}
\lastpage{384}

\makeatletter

\renewcommand{\hat}{\widehat}
\renewcommand{\widetilde}{\tilde}

\newtheorem{lemma}{Lemma}[section]
\newtheorem{prop}{Proposition}[section]
\newtheorem{definition}{Definition}[section]

\newtheorem{Cor}{Corollary}[section]

\newcommand{\Khat}{\widehat{\K}}

\newcommand{\indep}{\perp\hspace*{-6.2pt}\perp}

\newcommand{\real}[1]{{\mathbb R}^{#1}}
\newcommand{\spn}{\operatorname{span}}
\newcommand{\diag}{\operatorname{diag}}
\newcommand{\E}{\mathrm{E}}
\newcommand{\var}{\operatorname{var}}
\newcommand{\cov}{\operatorname{cov}}
\newcommand{\rank}{\operatorname{rank}}
\newcommand{\pred}{\mathrm{pred}}
\newcommand{\tr}{\operatorname{tr}}
\newcommand{\vecc}{\operatorname{vec}}
\newcommand{\I}{\mathbf I}

\newcommand{\rms}{\mathrm{RMS}}

\newcommand{\X}{{\mathbf X}}
\newcommand{\abf}{\mathbf a}
\newcommand{\f}{{\mathbf f}}
\newcommand{\bx}{{\mathbf x}}
\newcommand{\Pb}{{\mathbf P}}
\newcommand{\Q}{{\mathbf Q}}

\newcommand{\Xbar}{\bar{\X}}
\newcommand{\Xbb}{\mathbb{X}}

\newcommand{\D}{\mathbf D}
\newcommand{\K}{\mathbf K}

\newcommand{\ghat}{\hat{g}}

\newcommand{\Pbf}{{\mathbf P}}
\newcommand{\Qbf}{{\mathbf Q}}

\newcommand{\T}{\mathbf T}
\newcommand{\V}{{\mathbf V}}
\newcommand{\Vhat}{\widehat\V}

\newcommand{\W}{{\mathbf W}}
\newcommand{\Wi}{{\mathbf W}}
\newcommand{\What}{\widehat{\W}}
\newcommand{\Whati}{\widehat{\W}}
\newcommand{\Yhat}{\widehat{Y}}
\newcommand{\M}{{\mathbf M}}
\newcommand{\tn}{\widetilde{n}}
\newcommand{\Sbf}{{\mathbf S}}
\newcommand{\bbf}{{\mathbf b}}
\newcommand{\bbfhat}{\hat{\mathbf b}}
\newcommand{\A}{{\mathbf A}}
\newcommand{\B}{{\mathbf B}}
\newcommand{\G}{{\mathbf G}}
\newcommand{\R}{{\mathbf R}}
\newcommand{\Rhat}{\widehat{\R}}
\newcommand{\Rhatspice}{\Rhat_{\mathrm{spice}}}
\newcommand{\Rw}{\widehat{\R}_{\W}}
\newcommand{\Rwhat}{\widehat{\R}_{\What}}
\newcommand{\Ri}{\widehat{\R}_{\I}}
\newcommand{\Rg}{\widehat{\R}_{\widehat{\Deltabfs}}}
\newcommand{\Rd}{\widehat{\R}_{\mathrm{diag}}}
\newcommand{\J}{\mathbf J}
\newcommand{\Obf}{\mathbf O}

\newcommand{\etabf}{\bolds{\eta}}
\newcommand{\betabf}{\bolds{\beta}}
\newcommand{\epsilonbf}{\bolds{\epsilon}}
\newcommand{\varepsilonbf}{\bolds{\varepsilon}}

\newcommand{\Gammabf}{\bolds{\Gamma}}
\newcommand{\Gammabfshat}{\widehat{\bolds{\Gamma}}}
\newcommand{\Gammabfhat}{\widehat{\bolds{\Gamma}}}
\newcommand{\Deltabf}{\bolds{\Delta}}
\newcommand{\Thetabf}{\bolds{\Theta}}
\newcommand{\Deltabfhat}{\widehat{\bolds{\Delta}}}
\newcommand{\Deltabfs}{{\bolds{\Delta}}}
\newcommand{\Deltabfshat}{{\widehat{\bolds{\Delta}}}}
\newcommand{\Omegabf}{\bolds{\Omega}}
\newcommand{\Omegabfs}{{\bolds{\Omega}}}
\newcommand{\Omegabfhat}{\widehat{\bolds{\Omega}}}
\newcommand{\omegaS}{\omega_{\mathrm{spice}}}
\newcommand{\omegaD}{\omega_{\mathrm{diag}}}

\newcommand{\mubf}{\bolds{\mu}}
\newcommand{\mubfhat}{\hat{\mubf}}
\newcommand{\Fbb}{\mathbb{F}}
\newcommand{\rhobf}{\bolds{\rho}}
\newcommand{\rhobar}{\bar{\rho}}
\newcommand{\Xibf}{\bolds{\Xi}}
\newcommand{\Xibft}{\tilde{\bolds{\Xi}}}
\newcommand{\Bhat}{\widehat{\B}}
\newcommand{\Phibf}{\bolds{\Phi}}

\newcommand{\nubf}{\bolds{\nu}}

\newcommand{\syx}{\mathcal{S}_{Y|\X}}
\newcommand{\sexy}{\mathcal{S}_{E(\X|Y)}}

\newcommand{\spc}{{\mathcal S}}
\newcommand{\Zbb}{{\mathbb Z}}
\newcommand{\Fdmh}{\Phibf_{n}^{-1/2}}
\newcommand{\Fd}{\Phibf_{n}}
\newcommand{\Es}{\E_{\mathrm{sim}}}

\makeatother

\begin{document}
\begin{frontmatter}

\title{Estimating sufficient reductions of the predictors in abundant high-dimensional regressions}
\runtitle{Abundant sufficient dimension reduction}

\begin{aug}
\author[A]{\fnms{R. Dennis} \snm{Cook}\corref{}\thanksref{t1}\ead[label=e1]{dennis@stat.umn.edu}},
\author[B]{\fnms{Liliana} \snm{Forzani}\ead[label=e2]{liliana.forzani@gmail.com}}
\and
\author[A]{\fnms{Adam J.} \snm{Rothman}\thanksref{t2}\ead[label=e3]{rothman@stat.umn.edu}}
\runauthor{R. D. Cook, L. Forzani and A. J. Rothman}
\affiliation{University of Minnesota, Instituto de
Matem\'atica Aplicada del Litoral and~University of Minnesota}
\address[A]{R. D. Cook\\
A. J. Rothman\\
School of Statistics\\
University of Minnesota\\
Minneapolis, Minnesota 55455\\
USA\\
\printead{e1}\\
\hphantom{E-mail: }\printead*{e3}}
\address[B]{L. Forzani\\
Instituto de Matem\'atica Aplicada\\
\quad del Litoral \\
Facultad de Ingenier\'\i a~Qu\'\i mica\\
CONICET and UNL\\
G\"uemes 3450, (3000) Santa Fe\\
Argentina \\
\printead{e2}} 
\end{aug}

\thankstext{t1}{Supported by NSF Grants DMS-10-07547
and DMS-11-05650.}

\thankstext{t2}{Supported by NSF Grant DMS-11-05650.}

\received{\smonth{2} \syear{2011}}
\revised{\smonth{8} \syear{2011}}

%
\begin{abstract}
We study the asymptotic behavior of a~class of methods for sufficient
dimension reduction in high-dimension regressions, as the sample size
and number of predictors grow in various alignments. It is demonstrated
that these methods are consistent in a~variety of settings,
particularly in abundant regressions where most predictors contribute
some information on the response, and oracle rates are possible.
Simulation results are presented to support the theoretical conclusion.
\end{abstract}

%

\end{frontmatter}

\section{Introduction}\label{secintro}

There are many facets to the analysis of data in high dimensions, depending
on the type of application, some relying on dimension reduction, others
relying on variable selection and a~few employing both tactics. There has
been considerable interest in dimension-reduction methods for the regression
of a~real response $Y$ on a~random vector of predictors $\X\in\real{p}$
since the introduction of sliced inverse regression [SIR;
\citet{Li91}]
and sliced average variance estimation [SAVE; \citet{CooWei91}].
A common goal of these and many other methods is to reduce the dimension
of the predictor vector without loss of information about the response.
The aim is to estimate a~reduction $\R\dvtx\real{p}\rightarrow\real
{d}$, $d
\leq p$, with the property that $Y \indep\X|\R(\X)$ or,
equivalently, $\X|(Y,
\R(\X)) \sim\X|\R(\X)$ [\citet{Coo07}]. In this way $\R$ is
sufficient because
it
captures all the information about $Y$ that is available from $\X$. Sufficient
reductions are not determined uniquely by this definition because any bijective
transformation of $\R$ is also sufficient.

Nearly all methods for sufficient dimension reduction (SDR) restrict
attention to the class of linear reductions, which arise naturally in many
contexts. Linear reduction can be represented conveniently in terms of the
projection $\Pbf_{\spc}\X$ of $\X$ onto a~subspace $\spc\subseteq
\real{p}$. If $Y \indep\X|\Pbf_{\spc}\X$, then $\spc$ is called a~dimension-reduction subspace. Under mild conditions the intersection of any
two dimension-reduction subspaces is again a~dimension-reduction subspace
and that being so the central subspace $\syx$, defined as the intersection
of all dimension-reduction subspaces, is taken as the inferential target
[Cook (\citeyear{Coo94}, \citeyear{Coo98})]. A minimal
sufficient linear reduction is then of
the form
$\R(\X)=\etabf^{T}\X$, where $\etabf$ is any basis for $\syx$.

SDR has a~long history of successful application and is still an active
research area. Recent novel SDR methods include likelihood-based
sufficient dimension reduction [\citet{CooFor09}], kernel
dimension reduction [\citet{FukBacJor09}], shrinkage inverse
regression estimation [\citet{BonLi09}], dimension reduction for
nonelliptically distributed predictors [\citet{LiDon09},
\citet{DonLi10}], cumulative slicing estimation
[\citet{ZhuZhuFen10}], dimension reduction for survival models
[\citet{XiaZhaXu10}] and dimension reduction for spatial point
processes [\citet{GuaWan10}].
This body of work reflects three different but related frontiers in SDR:
extensions that require progressively fewer assumptions, development of
likelihood-based methods and
adaptations for specific areas of application. Almost all SDR
methods rely on traditional asymptotic reasoning for support, letting
the sample
size $n
\rightarrow\infty$ with $p$ fixed. They nearly all require the
inverse of
a~$p \times p$ sample covariance matrix and thus application is problematic
when $n < p$. Since accurate estimation of a~general $p \times p$
covariance matrix can require $n \gg p$ observations, it has seemed inevitable
that SDR methods would encounter estimation problems when $n$ is not
sufficiently large. \citet{ChiMar02}, \citet{LiLi04}
and others circumvented these issues by performing reduction in two
stages, first replacing the $p$ predictors with $p^{*} \ll n$ principal
components and then applying an SDR method to the regression of the
response on
the selected $p^{*}$ components. However, recent results on the
eigenvectors of sample covariance matrices in high-dimensional settings raise
questions on the value of such two-stage methods [see, e.g.,
\citet{JohLu09}]. \citet{CooLiChi07} proposed an SDR method
that avoids
computation of inverses and
reduces to partial least squares in a~special case. \citet{ChuKel10}
showed recently that the partial least squares estimator of the
coefficient vector in the linear regression of $Y$ on $\X$ is inconsistent
unless $p/n \rightarrow0$ and this raises questions about the behavior of
the Cook et al. SDR estimator when $n$ is not large relative to $p$.
\citet{LiYin08} used the least squares formulation of sliced inverse
regression originated by \citet{Coo04} to develop a~regularized
version that
allows $n < p$ and achieves simultaneous predictor selection and dimension
reduction. This seems to be a~promising method, but its asymptotic properties
are unknown and it may not work well when the regression is not sparse.
\citet{WuLi11} studied the asymptotic\vadjust{\goodbreak} properties of a~family of
SDR estimators,
using a~SCAD-type penalty for variable selection in sparse regressions
where the
number of relevant variables is fixed as $p \rightarrow\infty$. Their method
requires that $p/n \rightarrow0$ for consistency.

While sparsity is an important concept in high-dimensional regression,
not all
high-dimensional regressions are sparse.
For example, near-infrared reflectance is often measured at many wavelengths
to predict the composition of matter, like the protein content of a~grain.
There is not normally an expectation that only a~few wavelengths are
needed to
predict content. While some wavelengths may be better predictors than others,
it is the cumulative information provided by many wavelengths that is often
relevant. Partial least squares has been the dimension-reduction method of
choice in this type of regression. The regressions implied by this and other
nonsparse applications share similar characteristics: (1) the predictor
vectors are high-dimensional and typical area-specific analyses have employed
some type of dimension reduction; (2) while assessing the relative importance
of the predictors may be of interest, prediction is the ultimate goal; (3)
information on the response is thought to accumulate, albeit perhaps
slowly, as
predictors are added, and (4) sparsity is not a~driving notion.

In this article we introduce a~family of SDR methods for studying
high-dimensional regressions that
differs from past approaches in at least three important
ways. First, we do not require sparsity but rather we emphasize
\textit{abundant
regressions} where most of the predictors contribute some information
about the
response. In the logic of \citet{autokey23}, the bet-on-sparsity
principle arose because, to continue the metaphor, there is otherwise little
chance of a~reasonable payoff. We show in contrast that reasonable
payoffs can
be obtained in abundant regressions with prediction as the ultimate goal,
leading to a~contrasting \textit{bet-on-abundance} principle. Second,
SDR studies
have largely focused on properties of estimators of $\syx$. We bypass
this step
and instead consider the limiting behavior of estimators of the sufficient
reduction $\R(\X)$
itself, assuming that the dimension $d$ of $\R$ is fixed. More specifically,
letting $\Rhat$ denote an estimated reduction, we establish rates of
convergence
in the following sense. Let $\X_{N}$ denote a~new observation on~$\X$. If
$\Rhat(\X_{N}) - \R(\X_{N})= O_{p}(r(n,p))$ and if $r(n,p)
\rightarrow0$ as
$n,p \rightarrow
\infty$, then $\Rhat(\X_{N})$ is consistent for $\R(\X_{N})$ and
its convergence
rate is at least $r^{-1}$. Third, we integrate recent work on the
estimation of
high-dimensional covariance matrices into our approach. In particular, we
estimate a~critical matrix of weights by using sparse permutation invariant
covariance estimation (SPICE) as developed by \citet{Rotetal08}.

In sum, by considering the reduction $\R$ itself rather than the central
subspace $\syx$, we both introduce a~novel viewpoint for addressing dimension
reduction and develop theoretically grounded SDR methodology for
$n < p$ regressions where other methods either have no asymptotic
support or
must necessarily fail.\vadjust{\goodbreak}

We describe the model for our study and its sufficient reduction in
Section~\ref{secmodelstate}. The class of estimators that we use is
described in Section~\ref{secestimation}, and stabilizing
restrictions are
presented in Section~\ref{secsignal}. Sections~\ref{secWhat} and
\ref{secnormalerrors} contain theoretical conclusions for selected
estimators from the class described in Section~\ref{secestimation}.
Simulation results are presented in Sections~\ref{secsimW} and
\ref{secsimulations}. We turn to a~spectroscopy application in
Section~\ref{secdata} and a~concluding discussion is given in
Section~\ref{secdiscussion}. All proofs and additional simulation results
are available in a~supplemental
article [\citet{CooForRot}].

The following notational conventions will be used in our exposition. We
use $\real{p \times q}$ to denote the collection of all real $p \times
q$ matrices. We use $\|\A\|$ and $\|\A\|_{F}$ to denote the spectral
and Frobenius norms of $\A$. The largest and smallest eigenvalues of
$\A\in\real{p \times p}$ are denoted $\varphi_{\max}(\A)$ and
$\varphi_{\min}(\A)$. If $\A\in\real{p \times p}$, then $\diag(\A)
\in\real{p \times p}$ is the diagonal matrix with diagonal elements
equal to those of $\A$. $\vecc(\A)$ is the operator that maps
$\A\in\real{p \times q}$ to~$\real{pq}$ by stacking its columns. If
$\B\in\real{p \times q}$ and $\A\in\real{p \times p}$ is symmetric and
positive definite, then the operator that projects in the $\A$ inner
product onto\vspace*{1pt} $\spn(\B)$, the subspace spanned by the
columns of $\B$, has the matrix representation $\Pb_{\B(\A)} =
\B(\B^{T}\A\B)^{-1}\B^{T}\A$, and $\Q_{\B(\A)} = \I_{p} - \Pb_{\B(\A)
}$. $\Pb_{\B}$ indicates the projection onto $\spn(\B)$ in the usual
inner product. A basis matrix for a~subspace $\spc\subseteq\real{p}$ of
dimension $d$ is any matrix $\B\in\real{p \times d}$ whose columns
form a~basis for $\spc$. For nonstochastic sequences $\{a_{n}\}$ and~$\{b_{n}\}$,
we write $a_{n} \asymp b_{n}$ if there are constants $m$,
$M$ and $N$ such that $0 < m< |a_{n}/b_{n}|< M < \infty$ for all $n >
N$. Similarly, for stochastic sequences $\{a_{n}\}$ and~$\{ b_{n}\}$,
we write $a_{n} \asymp_{p} b_{n}$ if $a_{n}=O_{p}(b_{n})$ and
$b_{n}=O_{p}(a_{n})$. For $\A, \B\in\real{p \times p}$, $\A> \B$ means
that $\A-\B$ is positive definite. $\otimes$ denotes the Kronecker
product, and $X \sim Y$ means $X$ and $Y$ are equal in distribution.
Deviating slightly from convention, we do not index quantities by $n$
and $p$, preferring instead to avoid notation proliferation by giving
reminders from time to time.

\section{Model}\label{secmodelstate}

We assume throughout that the
data consist of independent observations $(Y_{i}, \X_{i})$,
$i=1,\ldots,n$, and that $p$ is increased by taking additional
measurements on
each of $n$
sampling units. We treat the process of
adding predictors as conditional on the response and so our
approach is based on
inverse reduction, $\X|(Y, \R(\X)) \sim\X|\R(\X)$, which is a~standard SDR
paradigm.
Limits as $n,p \rightarrow\infty$ are thus conditional on the
responses unless
specifically indicated otherwise.

\subsection{Inverse regression model}\label{secmodel}
The model that we employ is engendered by standard SDR assumptions. We first
review those assumptions briefly and then give a~statement of the model.

Classical methods like SIR require two
predominant and well-known conditions for estimation of $\syx$. The first,
called the \textit{linearity condition}, insists that $\E(\X|\etabf^T\X)$
be a~linear function of $\etabf^{T}\X$, where $\etabf$ is a~basis matrix
for~$\syx$.\vadjust{\goodbreak}
It must be valid only at a~true basis and not for all $\etabf\in\real{p \times d}$.
This condition is generally regarded as mild since it
holds to a~good approximation when $p$ is large [\citet{HalLi93}].
Let $\sexy$ denote the subspace spanned by $\E(\X|Y=y) -
\E(\X)$ as $y$ varies in the sample space of $Y$. Then the linearity
condition implies that $\sexy\subseteq\var(\X)\syx$, which is used
as a~basis for estimating a~subspace of $\syx$. The
second, called the \textit{coverage condition}, requires that $\sexy=
\var(\X)\syx$ and it
too is generally regarded as mild in many applications. While the
linearity and
coverage conditions are standard requirements for root-$n$ consistent methods
based on the inverse mean
function $\E(\X|Y)$, the actual performance of those methods depends
also on the
inverse variance function $\var(\X|Y)$. For instance,
\citet{BurCoo01} concluded based on simulations that SIR
works best when $\var(\X|Y)$ is nonstochastic,
and \citet{CooNi05} demonstrated analytically that SIR can be quite
inefficient when $\var(\X|Y)$ varies. We refer to the
requirement that $\var(\X|Y)$ be nonstochastic as the \textit{covariance
condition}. If the linearity and covariance conditions hold, then $\E
(\var(\X|\etabf^{T}\X)|Y)$ must be nonstochastic as well. This is
related to the usual covariance condition---$\var(\X|\etabf^{T}\X)$
is constant---used by SAVE and other second-order methods. See
\citet{Li91}, \citet{CooNi05} and \citet{LiDon09} for
further discussion of these
conditions.

We assume the linearity, coverage and covariance conditions as a~basis for our study. Under these conditions it can be shown
straightforwardly that $\var(\X)\syx= \Deltabf\syx$, where
$\Deltabf=
\var(\X|Y)$. Consequently, letting $\tilde{\Gammabf} \in\real{p
\times
d}$ be a~basis matrix for $ \Deltabf\syx$, we are led to the model
%
%
\begin{equation} \label{origmodel}
\X_{i} = \mubf+ \tilde{\Gammabf} \tilde{\bolds{\xi}}_{i}+
\varepsilonbf_{i},\qquad
i=1,\ldots,n,
\end{equation}
where\vspace*{1pt} $\mubf= \E(\X)$, the error vectors
$\varepsilonbf_{i}$ are independent copies of a~random vector
$\varepsilonbf\in\real{p}$ with mean 0 and variance
$\var(\varepsilonbf) = \Deltabf$, and $\tilde{\bolds{\xi}}_{i} =
\tilde{\bolds{\xi}}(Y_{i})$ is the $i$th instance of an unknown
vector-valued function $\tilde{\bolds{\xi}}(Y) \in\real{d}$ with
$\E(\tilde{\bolds{\xi}}) = 0$ that gives\vspace*{1pt} the coordinates
of $\E(\X|Y) - \E(\X)$ in terms of $\tilde{\Gammabf}$ and is
independent of $\varepsilonbf$. In this representation $\syx=
\Deltabf^{-1}\spn(\tilde{\Gammabf})$.

Neither $\tilde{\Gammabf}$ nor $\tilde{\bolds{\xi}}$ is
identified in
model~(\ref{origmodel}),
because for any conforming nonsingular matrix $\A$,
$\tilde{\Gammabf}\tilde{\bolds{\xi}}= (\tilde{\Gammabf}\A
)(\A^{-1}
\tilde{\bolds{\xi}})$,
leading to a~different parameterization. This nonuniqueness has been
mitigated
in past studies with $p$ fixed by requiring
$\tilde{\Gammabf}{}^{T}\tilde{\Gammabf} = \I_{d}$. However, that
parameterization is not workable when allowing $p \rightarrow\infty$,
and for this reason we adopt different restrictions.

The next step is to reparameterize model~(\ref{origmodel}) to satisfy
constraints that
will facilitate our development.
Specifically, we construct an equivalent
model
%
%
\begin{equation}\label{vecmodel}
\X_{i} = \mubf+ \Gammabf\bolds{\xi}_{i} + \varepsilonbf_{i},
\qquad i=1,\ldots,n,
\end{equation}
where $\spn(\Gammabf) = \spn(\tilde{\Gammabf})$, $\bolds{\xi}_{i} =
\bolds{\xi}(Y_{i})$ is the coordinate vector relative to the new basis
matrix $\Gammabf$, and we center the $\bolds{\xi}_{i}$'s in the sample
so that $\bar{\bolds{\xi}} = 0$. Let the rows\vadjust{\goodbreak} of $\Xibft\in\real{n
\times d}$ and $\Xibf\in\real{n \times d}$ be
$\tilde{\bolds{\xi}}{}^{T}_{i}$ and $\bolds{\xi}_{i}^{T}$,
$i=1,\ldots,n$. Without loss of generality, we impose on model
(\ref{vecmodel}) the constraints that (1) $n^{-1}\Xibf^{T}\M_{n}
\Xibf=\I_d$, where $\M_{n}$ is a~\textit{scaling matrix} that is
defined in Section~\ref{scaling}, and (2) $\Gammabf^{T}\W\Gammabf$ is a~diagonal matrix where $\W\geq0$ is a~symmetric $p \times p$ population
\textit{weight matrix} that is discussed in Section
\ref{secestimation}. To see how this is done starting from model
(\ref{origmodel}), let ${\mathbf T} = (n^{-1}\Xibft{}^{T} \M_n
\Xibft)^{-1/2}$ and\vspace*{1pt} let $\Obf\in\real{d \times d}$ be an
orthogonal
matrix so that $\Obf^{T} \T^{-1}\tilde{\Gammabf}{}^T \Wi\tilde{\Gammabf}
\T^{-1}\Obf$ is a~diagonal matrix. Then $\Xibft\tilde{\Gammabf}{}^{T} =
\Xibft\T\Obf\Obf^{T}\T^{-1}\tilde{\Gammabf}{}^{T}=\Xibf\Gammabf^{T}$,
where $\Xibf= \Xibft\T\Obf$ and $\Gammabf^{T} =
\Obf^{T}\T^{-1}\tilde{\Gammabf}{}^{T}$ satisfy the constraints by
construction. Model~(\ref{vecmodel}) can be represented also in matrix
form as
%
%
\begin{equation}\label{matmodel}
\Xbb= \mathbf1_{n}\mubf^{T} + \Xibf\Gammabf^{T} + {\mathbf e},
\end{equation}
where $\Xbb\in\real{n \times p}$ and ${\mathbf e}\in\real{n \times
p}$ have rows
$\X_{i}^{T}$ and $\varepsilonbf_{i}^{T}$, and ${\mathbf e}$ has mean
$0$ and variance
$\var(\vecc({\mathbf e}^{T})) = \I_{n}\otimes\Deltabf$.

Since bijective transformations of
a~sufficient reduction are themselves sufficient, we define $\R$ to be the
coordinates of the projection of $\X- \mubf$ onto $\spn(\Gammabf)$
in the
$\Deltabf^{-1}$ inner product:
%
%
\begin{equation}\label{popred}
\R(\X) =
(\Gammabf^{T}\Deltabf^{-1}\Gammabf)^{-1}\Gammabf^{T}\Deltabf
^{-1}(\X-
\mubf),
\end{equation}
where the first factor $(\Gammabf^{T}\Deltabf^{-1}\Gammabf)^{-1}$ stabilizes
$\E(\R)$ as $p \rightarrow\infty$.

\section{Estimation}\label{secestimation}

Without further structure it is not possible to use mod\-el~(\ref{vecmodel})
to estimate $\R(\X)$ since the coordinate vectors $\bolds{\xi
}_{i}$ are unknown.
However, estimation is possible by approximating the coordinate vectors as
$\bolds{\xi}_{i} \approx\bbf\f_{i}$, where $\bbf\in\real{d
\times r}$, $d
\leq r$, and $\f_{i} = \f(Y_{i})$ is the $i$th realization of a~known
user-selected vector-valued
function $\f(Y) \in\real{r}$.
Without loss of generality we center the sample $\sum_{i=1}^{n} \f
_{i}= 0$.
Let $\Fbb
\in\real{n \times r}$ be the matrix with rows $\f_{i}^{T}$, and
assume that
$\Phibf_{n} = \Fbb^T \Fbb/n > 0$ and $\Phibf= \lim_n \Fd> 0$,
unless $r=0$
and then of course $\Fbb$ is nil. We next describe the class of estimators
we use, postponing discussion of possible choices for $\f$ until
Section~\ref{choose-f}.

\subsection{Estimators}\label{estimators} The class of
reduction estimators $\Rhat$ that we study is based on using
estimators of $(\mubf, \bbf, \Gammabf)$ from a~subclass of the
family of
inverse regression estimators proposed by \citet{CooNi05}:
\[
(\mubfhat, \bbfhat, \Gammabfhat) = \arg\min\tr\{(\Xbb- \mathbf1
_{n}\mubf^{T} -
\Fbb\bbf^{T} \Gammabf^{T})\Whati(\Xbb- \mathbf1_{n}\mubf^{T} - \Fbb
\bbf^{T}\Gammabf^{T})^{T}\},
\]
where\vspace*{1pt} the minimization is over $\mubf\in\real{p}$,
$\Gammabf\in
\real{p
\times d}$ and $\bbf\in\real{d \times r}$ subject to the constraints that
$\bbf\Fd\bbf^{T}= \I_{d}$ and that $\Gammabf^{T}\What\Gammabf$
is a~diagonal matrix. A~particular estimator is determined by the choice of the
sample weight matrix $\What\in\real{p \times p}$ with $\W$ being a~population
version.
Some instances of~$\What$ that we introduce later correspond to $\W=
\I_{p}$,
$\Deltabf^{-1}$ and $ \diag^{-1}(\Deltabf)$. We next report the estimators
from this family. A sketch of the derivation is available in the supplemental
article [\citet{CooForRot}].\vadjust{\goodbreak}

It is easily seen that $\mubfhat= \Xbar$. Let $\Zbb= \Xbb-
\mathbf1_{n}\mubfhat^{T}$, and let $\Bhat= \Zbb^{T}\Fbb(\Fbb^{T}\Fbb
)^{-1} \in\real{p \times r}$ denote the matrix of regression
coefficients from the least\vspace*{1pt} squares fit of $\X$ on $\f$,
assuming that $n > r+1$. Also, let the columns of $\Vhat_{d} \in
\real{r \times d}$ be the first $d$ eigenvectors of
%
%
\begin{equation}\label{Khat}
\Khat= \Fd^{1/2}\Bhat^{T}\Whati\Bhat\Fd^{1/2} \in\real{r \times r},
\end{equation}
assuming that the $d$th eigenvalue of $\Khat$ is strictly larger
than the $(d+1)$st eigenvalue. This assumption will be true
with probability 1 for the
choices of~$\What$ considered here.
Then the estimators are $\bbfhat=
\Vhat_{d}^{T}\Fdmh$,
$\Gammabfhat= \Bhat\Fd^{1/2}\Vhat_{d}$,
$\Gammabfhat{}^{T}\Whati\Gammabfhat= \Vhat_{d}^{T}\Khat\Vhat_{d}$
and
%
%
\begin{equation} \label{samplered}
\Rhat(\X) = (\Gammabfhat{}^{T}\Whati\Gammabfhat)^{-1}
\Gammabfhat{}^{T}\Whati(\X- \Xbar).
\end{equation}
The diagonal elements of the diagonal matrix
$\Gammabfhat{}^{T}\Whati\Gammabfhat$ consist of the first (largest) $d$
eigenvalues of $\Khat$, and $\bbfhat$ is the coefficient matrix from
the OLS fit of $\Rhat(\X_{i})$ on~$\f_{i}$.

\subsection{Choice of $\f$}\label{choose-f}

Clearly, we should strive to choose an $\f(Y)$ so that, for some $\bbf
\in\real{d \times r}$, $\bolds{\xi}_{i}$ is well approximated by
$\bbf\f_{i}$ for $i=1,\ldots,n$. This intuition is manifested in
various calculations via the requirement that
$\rank(n^{-1}\Xibf^{T}\Fbb) =d$ for all~$n$. As an extreme first
instance, suppose that $\Xibf^{T}\Fbb= 0$. Then $\Bhat=
\Zbb^{T}\Fbb(\Fbb^{T}\Fbb)^{-1} = (\Gammabf\Xibf^{T}\Fbb+ {\mathbf
e}^{T}\Fbb)(\Fbb^{T}\Fbb)^{-1} = {\mathbf e}^{T}(\Fbb^{T}\Fbb)^{-1}$
and consequently~$\Bhat$ provides no information on $\spn(\Gammabf)$.
Since $n^{-1}\Xibf^{T}\Fbb= \widehat{\cov}(\bolds{\xi}(Y),\f(Y))$, this
requirement is essentially that the true coordinate vector
$\bolds{\xi}$ is sufficiently correlated with its approximation, and it
is equivalent to the condition derived by \citet{CooFor08} under a~PFC model with normal errors and $p$ fixed. In the remainder of this
article we assume that, for all $n$,
%
%
\begin{equation}\label{rank}
\rank(n^{-1}\Xibf^{T}\Fbb) =d \quad\mbox{and}\quad n^{-1}\Xibf
^{T}\Fbb= O(1) \qquad\mbox{as } n \rightarrow\infty,
\end{equation}
where the order is nonstochastic because we condition on the observed
values of $Y$ in our formal asymptotic calculations.

There are many ways to choose an appropriate $\f$ in practice.
Under model~(\ref{vecmodel}), each coordinate $X_{j}$, $j=1,\ldots,p$,
of $\X$ follows a~univariate linear model with predictor vector $\f(Y)$.
When $Y$ is quantitative we can use inverse response plots
[\citet{Coo98}, Chapter 10] of $X_{j}$ versus $Y$,
$j=1,\ldots,p$, to gain information about suitable choices for $\f$. When
$p$ is too large for a~thorough graphical investigation, which is the
case in
the context of this article,
we can consider $\f$'s that contain a~reasonably flexible set of basis
functions, like polynomial terms in $Y$. In some regressions there
may be a~natural choice for $\f$. For example, suppose that
$Y$ is categorical, taking
values in one of $m$ categories $C_k$, $k=1,\ldots,m$. We can
then set $r=m-1$ and specify the $k$th coordinate of $\f$ to be
the indicator function $J(y \in C_k)$. Another option consists of
``slicing'' the observed range of a~continuous $Y$ into $m$ bins
(categories) $C_k$,
$k=1,\ldots,m$. A rudimentary set of basis
functions can then be constructed by specifying the $k$th coordinate of
$\f$ as for the case of a~categorical $Y$. This has the effect of
approximating each
conditional mean $\E(X_{j}|Y)$ as a~step
function of $Y$ with $m$ steps. This option is of particular
interest because it exactly reproduces SIRs estimate of $\syx$ in the
traditional large-$n$ setting [\citet{CooFor08}]. Cubic splines
with the endpoints of the bins as the knots are a~more responsive
option that is
less prone to loss of intra-slice information.

\section{Universal context}\label{secsignal}

Our goal is to study the limiting behavior of\break $\Rhat(\X_{N}) - \R(\X_{N})$,
where $\R$ and $\Rhat$ are given by~(\ref{popred}) and~(\ref{samplered}).
This difference
depends on the choice of $\What$ and the behavior of the true
reduction $\R$ as
$p \rightarrow\infty$. In this section we measure and constrain the
interaction between~$\W$ and the model. We also place weak constraints
on $\R$
to help ensure well-behaved limits. The context described in this
section will
be assumed to hold throughout this article, without necessarily being
declared in
formal statements. We will revisit these constraints occasionally
during the
discussion, particularly when some of them are implied by other structure.

\subsection{Scaling matrix $\M_{n}$} \label{scaling}

Since bijective transformations of sufficient reductions are also
sufficient, we need to have $\Rhat(\X_{N})$ and $\R(\X_{N})$ in
the same scale to ensure that their difference $\Rhat(\X_{N}) - \R
(\X_{N})$
is a~useful measure of agreement. This can be accomplished by choosing the
scaling matrix $\M_{n} = \Pbf_{\Fbb}$ so that the first model constraint
stated following~(\ref{vecmodel}) becomes $n^{-1}\Xibf^{T}\M
_{n}\Xibf=
n^{-1}\Xibf^{T}\Pbf_{\Fbb}\Xibf\in\real{d \times d}$,
which is a~rank-$d$ matrix by
(\ref{rank}). This choice ensures that
$\Gammabf$ is the lead term in the asymptotic expansion of
$\Gammabfhat$, which places
$\Rhat(\X_{N})$ and $\R(\X_{N})$ on the same scale.

\subsection{Signal rate}

We define the \textit{signal rate} to be a~positive monotonically increasing
function $h(p) = O(p)$ that measures the rate of increase in the population
signal: let $\G_{h} = \Gammabf^{T}\Wi\Gammabf/h(p)$, which is a~diagonal matrix
because $\Gammabf^{T}\Wi\Gammabf$ is diagonal by construction. Then $h(p)$
is chosen to meet the requirement
%
%
\begin{equation}\label{h}
\lim_{p \rightarrow\infty} \G_{h} = \G\in\real{d \times d},
\end{equation}
where $\G$ is a~diagonal matrix that is assumed to be positive
definite with
finite elements.
The signal rate is not needed for computation of $\Rhat$ but it will
play a~key
role in later developments. It depends via $\W$ on the specific estimator
selected, although this is not indicated notationally. When $\W> 0$ the
bounds $\varphi_{\min}(\W)\Gammabf^{T}\Gammabf\leq\Gammabf
^{T}\Wi\Gammabf\leq
\varphi_{\max}(\W)\Gammabf^{T}\Gammabf$ can be used to aid
intuition on the
magnitude of $h(p)$ by presuming properties of $\Gammabf$ and $\W$. For
example, consider regressions in which $\varphi_{\min}(\W)$ and
$\varphi_{\max}(\W)$ are bounded away from 0 and $\infty$ as $p
\rightarrow
\infty$ and a~positive fraction $g$, $0 < a~\leq g \leq1$, of the
rows of
$\Gammabf$ is sampled from a~multivariate density with finite second
moments and
the other rows of $\Gammabf$ are all zero vectors. Then the number of nonzero
rows $pg \asymp p$, $h(p) \asymp p$, and we say that the regression has an
\textit{abundant signal}. Regressions with \textit{near abundant signals}
like $h(p)
\asymp p^{2/3}$ may often perform similarly in practice to regressions with
abundant signals. It is technically possible to have regressions in which
$p = o(h(p))$, but we would not normally expect this in practice. On the
other extreme, sparse regressions are those in which $h(p) \asymp1$,
so only
finitely many predictors are relevant or the signal accumulates very
slowly as
$p \rightarrow\infty$. Regressions with \textit{near sparse signals}
have, say,
$h(p) = o(p^{1/3})$. Typically we will use $h(p)$ in definitions and formal
statements, but use the abbreviated form $h$ otherwise. We assume in
the rest
of this article that~(\ref{h}) holds.

\subsection{Limiting reduction}\label{seclimitred}
Our second
requirement is that, as $p \rightarrow\infty$, the spectral norm $\|
{\var}(\R)\|$
converges to a~finite constant, which may be 0. If this is not so, then the
variance of some linear combinations $\abf^{T}\R$ will diverge as $p
\rightarrow
\infty$, and the very notion of dimension reduction for large-$p$ regressions
becomes problematic.
%
%
\begin{definition} A reduction $\R(\X)$ is stable if ${\lim
_{p\rightarrow
\infty}}\|{\var}(\R(\X))\| < \infty$.
\end{definition}

The following lemma gives a~sufficient condition for stability that incorporates
the signal rate. In preparation, define
%
%
\begin{equation} \label{rho}
\rhobf= \Wi^{1/2}\Deltabf\Wi^{1/2} \in\real{p \times p},
\end{equation}
which is manifested in various asymptotic expansions as a~measure of the
agreement between $\W$ and $\Deltabf^{-1}$.
%
%
\begin{lemma}\label{lemmastable} If $\|\rhobf\| = O(h(p))$, then reduction
(\ref{popred}) is stable.
\end{lemma}

According to Lemma~\ref{lemmastable}, if $\W= \Deltabf^{-1}$, then the
reduction is stable regardless of $h$. All regressions in which
$\|\rhobf\|$ is bounded yield stable reductions. Bounded eigenvalues
have been required in various studies to avoid ill-conditioned covariance
matrices [see, e.g., \citet{BicLev08N1} and \citet
{Rotetal08}]. More generally, the requirement for a~stable reduction is that
$h(p)$ must increase at a~rate that is no less than the rate at which
$\|\rhobf\|$ increases. In particular, if $h = O(1)$, then the sufficient
condition of Lemma~\ref{lemmastable} requires that $\|\rhobf\|$ is bounded.

The rates that we develop depend on the functions $\kappa(n,p) =
[p/\{h(p)n\}]^{1/2}$ and $\psi(n,p,\rhobf) = \|\rhobf\|_{F}/\{
h(p)\sqrt{n}\}$.
The interpretation and roles of these functions\vadjust{\goodbreak} will be discussed
later; for now
we state across-the-board requirements that, as $n,p \rightarrow\infty$,
%
%
\begin{equation}\label{kappa}
\kappa(n,p) = O(1) \quad\mbox{and}\quad \psi(n,p,\rhobf) = O(1).
\end{equation}
These functions will nearly always be written without their arguments. We
assume in the rest of this article that all regressions are stable and that
the orders given in~(\ref{kappa}) hold.

\subsection{\texorpdfstring{Joint constraints on the weights $\W$ and errors $\varepsilon$}{Joint constraints on the weights W and errors epsilon}}
So far we have assumed only that the errors $\varepsilonbf_{i}$ are
independent copies of the random vector~$\varepsilonbf$ which has mean
0 and
variance $\Deltabf$. We impose two additional constraints to ensure
well-behaved limits: as $p \rightarrow\infty$,
\begin{longlist}[(W.2)]
\item[(W.1)] $\E(\varepsilonbf^{T}\W\varepsilonbf)
= O(p)$,
and
\item[(W.2)] $\var(\varepsilonbf^{T}\W\varepsilonbf
) =
O(p^{2})$.
\end{longlist}
Condition (W.1), which is equivalent to $\tr(\rhobf)/p = O(1)$, seems
quite mild.
For example, if we use unweighted fits with $\What= \W= \I_{p}$,
then this
condition is simply that the average error variance $\tr(\Deltabf)/p$ is
bounded.
Condition~(W.2) can be seen as placing constraints on $\W$ and on the fourth
moments of~$\varepsilonbf$. It too seems mild, although it is more constraining
than the first condition. The following lemma describes a~few settings
in which
condition (W.2) holds.
%
%
\begin{lemma}\label{lemmaW} Let $\epsilonbf= (\epsilon_{i}) =
\Deltabf^{-1/2}\varepsilonbf$ so that $\E(\epsilonbf) = 0$ and
$\var(\epsilonbf)
= \I_{p}$. Let $\phi_{ij} = \E(\epsilon_{i}^{2}\epsilon_{j}^{2})$,
$i,j=1,\ldots,p$. Then:
\begin{longlist}
\item If $\varepsilonbf\sim N(0,\Deltabf)$, then condition
\textup{(W.1)}
implies condition \textup{(W.2)}.
\item If $\W= \Deltabf^{-1}$ and if $\phi_{ij} \leq\phi<
\infty$
as $p \rightarrow\infty$, then condition \textup{(W.2)} holds.
\item If \textup{(a)} the elements $\epsilon_{i}$ of $\epsilonbf$ have
symmetric distributions or are independent and if \textup{(b)} $\phi_{ij} \leq
\phi<
\infty$ as $p \rightarrow\infty$, then condition \textup{(W.1)} implies
condition \textup{(W.2)}.
\end{longlist}
\end{lemma}

We assume in the rest of this article that conditions (W.1) and (W.2) hold.

\section{\texorpdfstring{$\What$ converging in spectral norm}{W converging in spectral norm}}\label{secWhat}

In this section we describe our results for the limiting reduction when
$\What$
converges in the spectral norm. Specifically, we require the following two
conditions:
\begin{longlist}[(S.2)]
\item[(S.1)] There exists a~population weight matrix $\W\in\real{p
\times
p}$ so that the spectral norm of $\Sbf\equiv\W^{-1/2}(\What- \W)\W^{-1/2}$
converges to 0 in probability at rate at most $\omega^{-1}(n,p)$ as
$n, p
\rightarrow\infty$; equivalently, $\|\Sbf\| = O_{p}(\omega(n,p))$
as $\omega
\rightarrow0$.
\item[(S.2)] $\|\E(\Sbf^{2})\| = O(\omega^{2}(n,p))$ as $\omega
\rightarrow
0$.\vadjust{\goodbreak}
\end{longlist}
These conditions are implied by the stronger condition that $\E(\|\Sbf
^{2}\|) =
O(\omega^{2})$, (S.1) following from Chebyshev's inequality and (S.2)
arising since
$\|\E(\Sbf^{2})\| \leq\E(\|\Sbf^{2}\|)$. All of the
weight matrices discussed in this article can satisfy these conditions,
as well
as other weight matrices that we have considered, so these conditions
do not
seem burdensome. For ease of reference we denote the corresponding estimator
as $\Rwhat$.

\subsection{Theoretical results}\label{secTRWhat}

We state and discuss one of our main results in this section. In preparation,
let
%
%
\begin{equation}\label{K}
\K= n^{-2}\Phibf_{n}^{-1/2}\Fbb^{T}\Xibf\G_h \Xibf^{T}\Fbb\Phibf
_{n}^{-1/2} \in\real{r
\times r},
\end{equation}
where $\G_{h}$ is as defined in~(\ref{h}), let $\rhobar=
\tr(\rhobf)/p$ and define the random vector $\nubf= \R_{\W
}(\varepsilonbf_{N})
- \R(\varepsilonbf_{N}) \in\real{d}$, where $\R_{\W
}(\varepsilonbf_{N}) =
(\Gammabf^{T}\W\Gammabf)^{-1}\Gammabf^{T}\W\varepsilonbf_{N}$ is
the population
reduction using $\W$ applied to the error $\varepsilonbf_{N}$ of a~new
observation and $\R(\varepsilonbf_{N})$ is the targeted population reduction
(\ref{popred}) applied to the same error. The vector $\nubf\in\real{d}$
measures the population-level agreement between the user-selected reduction
$\R_{\W}$ and the ideal reduction. It has mean \mbox{$\E(\nubf) = 0$} and variance
%
%
\begin{equation}\label{nu}
\var(\nubf) = \G_{h}^{-1/2}\biggl\{\frac{\bolds{\gamma
}^{T}\rhobf\bolds{\gamma}-
(\bolds{\gamma}^{T}\rhobf^{-1}\bolds{\gamma
})^{-1}}{h(p)}\biggr\}\G_{h}^{-1/2} \in\real{d
\times d},
\end{equation}
where the semi-orthogonal matrix $\bolds{\gamma}=
\Wi^{1/2}\Gammabf(\Gammabf^{T}\Wi\Gammabf)^{-1/2}$.

The next proposition describes asymptotic properties of $\Khat\in
\real{r
\times r}$ and~$\Rwhat$, where $\Khat$ was defined in~(\ref{Khat}).
%
%
\begin{prop}\label{propWhat1} \label{propWhat2} Assume conditions
\textup{(S.1)} and
\textup{(S.2)}.
Then:
\begin{longlist}
\item If $\Sbf=0$, then $h^{-1}(p)\E( \Khat)= \K+
\kappa^{2}\rhobar\I_r $.
\item $h^{-1}(p)(\Khat-\K) = O_{p}(\kappa) + O_{p}(\omega)$.
\item If, in addition, $\|\rhobf\| = O(h(p))$, then
\[
\Rwhat(\X_{N}) - \R(\X_{N}) = \nubf+ O_{p}(\kappa) + O_{p}(\psi)
+ O_{p}
(\omega),
\]
where $\kappa$ and $\psi$ were defined in~(\ref{kappa}).
\end{longlist}
\end{prop}

This proposition shows that the asymptotic properties of $\Khat$ and
$\Rwhat$ depend on four quantities, each with a~different role. We
defer discussion of the order $O_{p}(\omega)$ that measures the
asymptotic behavior of $\What$ to later sections. In the rest of this
section we concentrate on the remaining\vspace*{1pt} \mbox{terms---$\nubf$, $O_{p}(\kappa)$}
and $O_{p}(\psi)$---by assuming that $\What$ is
nonstochastic, so $\What= \W$, $\Sbf= 0$ and $\Rwhat= \Rw$. This
involves no loss since none of these terms depends on $\Sbf$.

Terms\vspace*{1pt} involving $\kappa$ affect both estimation $\Khat$ and prediction
$\Rw(\X_{N}) - \R(\X_{N})$. If $\kappa$ were to diverge, then, from
Proposition~\ref{propWhat1}(i), eventually $\E(\Khat)/h$ will look
like a~diagonal\vadjust{\goodbreak} matrix and there will be little signal left, which is why we
imposed the universal condition~(\ref{kappa}) that $\kappa= O(1)$. If
$\rhobar$ diverges, then again $\E(\Khat)/h$ will resemble a~diagonal matrix,
but
this is prohibited by universal condition (W.1). Here we also see the
role of
the rank condition~(\ref{rank}) introduced at the beginning of
Section~\ref{choose-f}.
If $\rank(n^{-1}\Xibf^{T}\Fbb) <d$, then $\rank(\K) < d$
and again some signal will lost. In the extreme, if $\Xibf^{T}\Fbb=
0$, then
$h^{-1}\E(\Khat) = \kappa^2\rhobar\I_r$ and all information on
$\spn(\Gammabf)$ is lost.

If $\kappa\rightarrow0$,
then $hn$ increases faster than $p$ increases. In effect there is a~synergy between the signal\vspace*{1pt} rate and the sample size. If the regression is
abundant, so $h \asymp p$, then $\kappa\asymp n^{-1/2}$. Useful
results can
also be obtained when $h \asymp p^{2/3}$, because then $\kappa\asymp
p^{1/6}n^{-1/2}$,
which may be small in some regressions. In sparse regressions where $h
\asymp1$,
$\kappa\asymp(p/n)^{1/2}$, and we are\vspace*{1pt} back to the usual requirement that
$p/n \rightarrow0$ for consistency of $\Khat$. Equally important, since
$\kappa$ does not depend on the weight matrix, the results indicate
that it
may
not be possible to develop from the family of estimators considered
rates of
convergence that are faster than $\kappa^{-1}$.

The $\nubf\in\real{d}$ and $O_{p}(\psi)$ terms arise from
prediction. The
first term $\nubf$ has a~characteristic that is different from the others
because it does not depend on $n$. Since $\var(\nubf) \leq
\G_{h}^{-1}\|\rhobf\|/h$, the sufficient stability requirement $\|
\rhobf\| =
O(h)$ of Lemma~\ref{lemmastable} guarantees that $\var(\nubf)$ is
bounded as $p
\rightarrow\infty$. While the contribution of $\nubf$ will be
negligible if
$\var(\nubf)$ is sufficiently small, for the best results we should have
$\var(\nubf) \rightarrow0$ as $ p \rightarrow\infty$. Let
$(\bolds{\gamma},
\bolds{\gamma}_{0})$ be an orthogonal matrix. Using a~result from
Cook and Forzani\vspace*{1pt}
[(\citeyear{CooFor09}), eq. (A.4)],
$
\bolds{\gamma}{}^{T}\rhobf\bolds{\gamma}- (\bolds
{\gamma}^{T}\rhobf^{-1}\bolds{\gamma})^{-1} =
\bolds{\gamma}^{T}\rhobf\bolds{\gamma}_{0}(\bolds
{\gamma}_{0}^{T}\rhobf\bolds{\gamma}_{0})^{-1}\bolds
{\gamma}_
{0}^{T}\rhobf\bolds{\gamma}_{0}.
$
Consequently, $\var(\nubf) = 0$ when $\spn(\bolds{\gamma})$ is
a~reducing subspace
of $\rhobf$, even if $h \asymp1$. This result is similar to Zyskind's
(\citeyear{Zys67})
classical findings about conditions for equality of the best and simple least
squares linear estimators in linear models. If $\W= \Deltabf^{-1}$,
then $\rhobf
= \I_{p}$ and $\spn(\bolds{\gamma})$ is trivially a~reducing
subspace of $\rhobf$. If
$\Deltabf$ is a~generalized inverse of $\W$, $\W\Deltabf\W= \W$,
and if
$\Gammabf^{T}\W\Gammabf> 0$, then again $\spn(\bolds{\gamma
})$ reduces $\rhobf$ and
$\var(\nubf) = 0$.

Turning to $O_{p}(\psi)$, since $\|\rhobf\|_{F} \leq\sqrt{p}\|
\rhobf\|$, it
follows that $\psi\leq\kappa
\|\rhobf\|/\sqrt{h}$. Consequently, if $\|\rhobf\|=O(\sqrt{h})$, then
$O_{p}(\kappa) + O_{p}(\psi) = O_{p}(\kappa)$ and $\Rw(\X_{N}) -
\R(\X_{N}) = \nubf+ O_{p}(\kappa)$. In a~worst-case scenario, if
$\psi= \kappa\|\rhobf\|/\sqrt{h}$ and $\|\rhobf\|/\sqrt{h}$ diverges,
then $O_{p}(\kappa) + O_{p}(\psi) = O_{p}(\sqrt{p}/\sqrt{n})$
because of the
requirement $\|\rhobf\| = O(h)$ in Proposition~\ref{propWhat1}(iii),
and these
terms reduce to the usual condition that $p/n \rightarrow0$.

In sparse regressions, the condition $\kappa\rightarrow0$ reduces to $p/n
\rightarrow0$ and $\|\rhobf\| = O(h)$ means that $\|\rhobf\|$ must be
bounded. These two conditions imply that \mbox{$\psi\rightarrow0$}.
Consequently, the best results in sparse regressions will be achieved when
(a) $\spn(\bolds{\gamma})$ is a~reducing subspace of $\rhobf$,
(b) $p/n
\rightarrow0$ and (c) $\|\rhobf\|$ is bounded.
In nonsparse settings $\Rw$ will yield the best results when $\|\rhobf
\|$ is
bounded and the regression is abundant or near abundant. We summarize some
implications of these conditions in the following corollary.\vadjust{\goodbreak}
%
%
\begin{Cor}\label{corW21} Assume that $\|\rhobf\| = O(1)$ and that
$\Sbf= 0$.
Then:
\begin{longlist}
\item $\Rw(\X_{N}) - \R(\X_{N}) = O_{p}(h^{-1/2}) + O_{p}(\kappa)$. If
the regression\vspace*{1pt} is abundant, then $\Rw(\X_{N}) - \R(\X_{N})
= O_{p}(n^{-1/2})$.
\item If $\spn(\bolds{\gamma})$ is a~reducing subspace
of $\rhobf$, then
$\Rw(\X_{N}) - \R(\X_{N}) = O_{p}(\kappa)$. If, in addition, the
regression is
abundant, then $\Rw(\X_{N}) - \R(\X_{N}) = O_{p}(n^{-1/2})$.
\end{longlist}
\end{Cor}

The conclusions of Proposition~\ref{propWhat2} and its Corollary~\ref{corW21}
hold as $n \rightarrow\infty$ and $p \rightarrow\infty$ in the required
relationships and as $n \rightarrow\infty$ with $p$ fixed. In the latter
case
the results simplify to those given in the next corollary.
%
%
\begin{Cor}\label{corW22} If $p = O(1)$ and $\Sbf= 0$, then:
\begin{longlist}
\item $\Rw(\X_{N}) - \R(\X_{N}) = \var(\nubf) + O_{p}(n^{-1/2})$.

\item If, in addition, $\spn(\bolds{\gamma})$ is a~reducing subspace of
$\rhobf$, then $\Rw(\X_{N}) - \R(\X_{N}) = O_{p}(n^{-1/2})$.
\end{longlist}
\end{Cor}

It may be clear from the previous discussion that the best possible
rate is
achieved when $\What= \W= \Deltabf^{-1}$ and then $\Rw(\X_{N}) -
\R(\X_{N}) =
O_{p}(\kappa)$. Consequently, we refer to $\kappa^{-1}$ as the
\textit{oracle rate}.

\subsection{Simulations}\label{secsimW}
We conducted a~simulation study with $\What= \W= \I_{p}$ to show the
importance of $\nubf$, to demonstrate that Proposition~\ref{propWhat2} and
Corollaries~\ref{corW21} and~\ref{corW22} give good qualitative
characterizations of the limiting behavior of $\Ri(\X_{N}) - \R(\X
_{N})$, and to
provide some intuition into the nonasymptotic behavior of $\Ri$. The simulated
data were generated using a~simple version of model~(\ref{vecmodel}):
$\X=
\Gammabf Y + \varepsilonbf$, $Y \sim N(0,1)$, $\Gammabf\in\real{p}$, was
constructed as a~vector of standard normal random variables and
$\Deltabf$ is a~diagonal matrix. Specific scenarios may differ on $n$, $p$ and the
choice of
$\Deltabf$. In any case, the true model then has $d=1$. We confined attention
to
regressions with $d=1$ because we found that there is nothing in
principle to
be
gained from settings with $d > 1$. For each sample size $n$ we used the
estimators
described at the end of Section~\ref{secestimation} with $\What= \I_{p}$
and $Y^{j}$, $j=1,2,3,4$, as the elements of $\f$, so $\rhobf=
\Deltabf$ and $r
= 4$.

For each data set constructed in this way, we determined $\Ri(\X
_{N,j})$\break and~$\R(\X_{N,j})$
at $j=1,\ldots,100$ new data points generated from the original model. We
summarized each data set by computing the absolute sample correlation between
$\Ri$ and $\R$ and the sample standard deviation of the difference
$\Ri- \R$
over the 100 new data points. This process was repeated at least 400
times for
$p \leq300$, 200 times for $300 < p < 800$ and 50 times for $p > 800$. The
average absolute correlation and average standard deviation were used
as overall
summary statistics.
Our decision to use a~diagonal matrix for~$\Deltabf$ was based on the
observation that the asymptotic results depend on~$\Deltabf$ largely via
$\|\rhobf\|$ and with $\What= \I_{p}$, $\|\rhobf\| =
\|\Deltabf\|$.

%
\begin{figure}
\begin{tabular}{@{}cc@{}}

\includegraphics{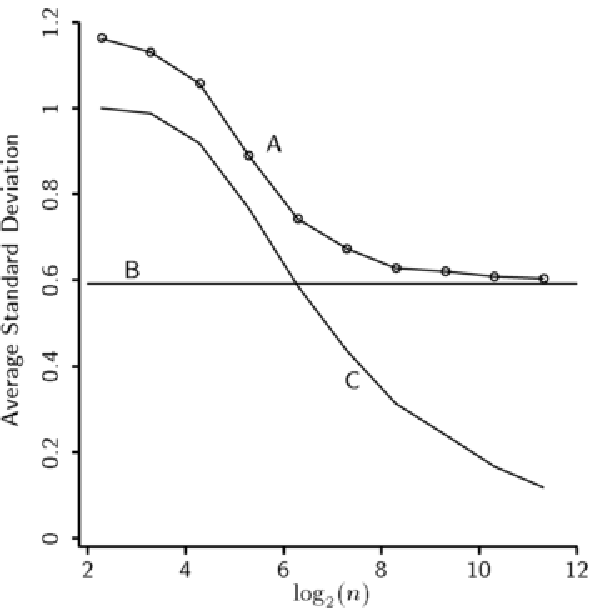}
 & \includegraphics{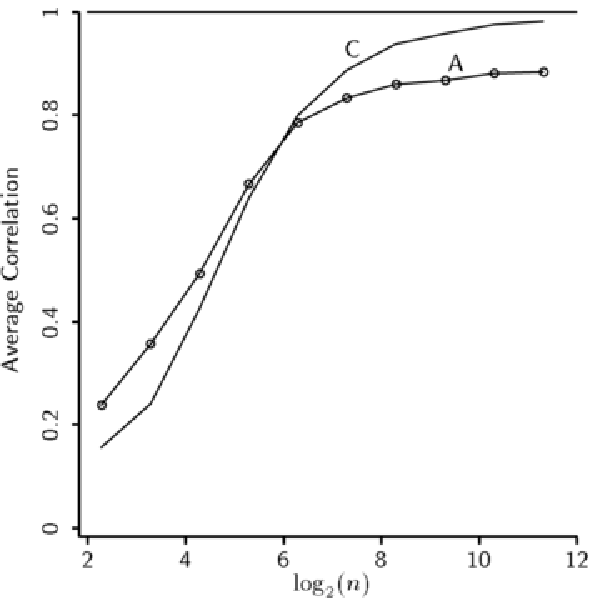}\\
(a) & (b)
\end{tabular}
\caption{Results of the first simulation described in Section
\protect\ref{secsimW} with fixed $p=50$, $\What=\I_{p}$, $Y \sim
N(0,1)$ and $\diag(\Deltabf) \sim\operatorname{Uniform}(1,101)$. A:
$\Gammabf$ generated as a~vector of $N(0,1)$ random variables. B:
Theoretical lower bound on the standard deviation of $\Ri- \R$ for
curve A. C:~$\Gammabf= (8.2,0,\ldots,0)^{T}$. \textup{(a)} Standard deviation of
$\Ri-\R$; \textup{(b)} correlation $(\Ri,\R)$.} \label{SDWI}
\end{figure}

\textit{Simulation} I: \textit{Illustrations of
Proposition}~\ref{propWhat2}.
In the first
set of simulations we fixed $p=50$ and generated $\Gammabf$ once at
the outset,
giving $\|\Gammabf\| = 8.2$. The diagonal elements of $\Deltabf$ were 50
regularly spaced values between~1 and~101. According to
Corollary~\ref{corW22}(i), as $n \rightarrow\infty$ the variance of
$\Ri- \R$
should converge to $\var(\nubf) = 0.595^2$, where the value was
determined by
substituting the simulation parameters into~(\ref{nu}). This
conclusion is
supported by the results in Figure~\ref{SDWI}(a), where we see that
the standard
deviation curve A converges nicely to the predicted asymptotic value
marked by
line segment~B. Curve~C in Figure~\ref{SDWI}(a) was generated by
setting $\Gammabf
= (8.2,0,\ldots,0)^{T}$, so $\spn(\bolds{\gamma})$ is now a~reducing subspace of
$\rhobf= \Deltabf$. In this case Corollary~\ref{corW22}(ii)
predicts that
$\Ri- \R$ will converge to 0 at the usual root-$n$ rate. That
conclusion is
supported by curve C. Curves A and C in Figure~\ref{SDWI}(b)
correspond to those
in Figure~\ref{SDWI}(a), except the average correlation is plotted on
the vertical
axis. Evidently the correlations for curve A are converging to a~value
close to
0.9, while those for curve C are converging to 1. Curve A suggests that when
$\var(\nubf) >0$ we might still gain useful results in practice if the
correlation is sufficiently large.

\textit{Simulation} II: \textit{Bounded and unbounded error variances.} The
structure of the
second set of simulations was like the first set, except we varied
$n=p/2$ and,
at the start of each replication, $\Gammabf$ and $Y$ were generated
anew and
the diagonal elements of $\Deltabf$ were sampled from the uniform distribution
on the interval $(1,u)$, where $u=101$ or $u=p+1$. The regeneration at the
start of each replication was used to avoid tying results to a~particular
parameter configuration. For this simulation the ratio of the systematic
variation in $\X$ to its total variation is $\var^{-1/2}(\X)\var(\E
(\X|\Gammabf,
Y, \Deltabf))\var^{-1/2}(\X) = (2/(u+3))\I_{p}$, where the moments
were computed
over the simulation distributions of $\X$, $Y$, $\Gammabf$ and
$\Deltabf$.
Consequently, it seems fair to characterize the signal as weak when $u
\ge101$,
since then the systematic variation accounts for less than 2\% of the total
variation.

\begin{figure}
\begin{tabular}{@{}cc@{}}

\includegraphics{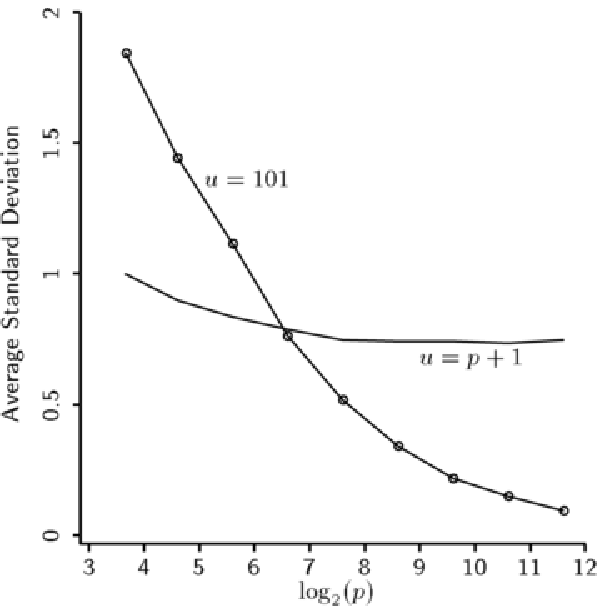}
 & \includegraphics{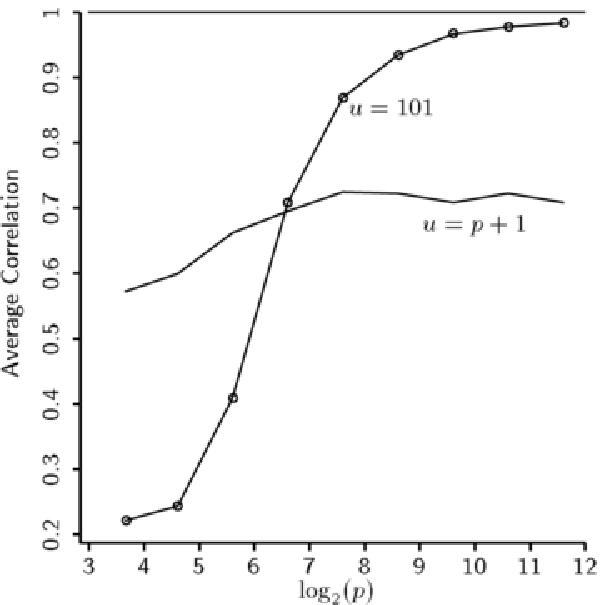}\\
(a) & (b)
\end{tabular}
\caption{Results from simulation II described in Section
\protect\ref{secsimW} with $n=p/2$, $\W=\I_{p}$, $Y \sim N(0,1)$,
$\Gammabf\sim N(0,\I_{p})$ and $\diag(\Deltabf) \sim
\operatorname{Uniform}(1,u)$. \textup{(a)} Standard deviation of $
\Ri-\R$; \textup{(b)} correlation $(\Ri,\R)$.} \label{SDU1p}
\end{figure}

The regression is abundant and $\kappa= o(1)$ in this simulation. Additionally,
when $u=101$, $\|\rhobf\|$ is bounded and the results of
Corollary~\ref{corW21}(i) apply, $\Ri- \R= O_{p}(n^{-1/2})$. The simulation
results for this case, which are shown by the curves labeled ``$u=101$''
in Figure~\ref{SDU1p}, support the claim that standard deviation of
%
%
$\Ri-\R$ is
converging to 0 and the correlation between $\Ri$ and $\R$ is
approaching 1. On
the other hand, when $u = p+1$, $\|\rhobf\|$ is unbounded and the
results of
Corollary~\ref{corW21} do not apply.
However, we still might expect the standard deviation and correlation to
converge to values away from~0 and~1. The simulation results shown by the
curves labeled $u=p+1$ in Figure~\ref{SDU1p} sustain this expectation.

\subsection{SPICE}\label{secspice}
Beginning with a~little background, we turn in this section to weight
matrices $\What\in\real{p \times p}$ selected by using SPICE.

Restricting attention to population weight matrices that are equal to
the inverse of the error covariance matrix, $\W= \Deltabf^{-1}$, allows
for the application of some modern regularized inverse covariance
estimators with reasonable convergence rates of $\What\in\real{p \times
p}$ as $n, p \rightarrow\infty$. For example, rates of convergence have
been established for banding or tapering [\citet{BicLev08N1}] and
element-wise thresholding [\citet{BicLev08N2}] of the sample
covariance matrix. However, using these approaches to estimate the
population weight matrix is problematic. Banding/tapering the sample
covariance matrix is not invariant under permutations of variable
labels, and in finite samples, both banding/tapering and thresholding
the sample covariance matrix may produce an estimator of the covariance
matrix that is not positive definite. To avoid these problems, we
further restrict our attention to covariance estimators that are
positive definite in finite samples and are invariant to permutations
of variable labels. Several authors have recently analyzed the
high-dimensional inverse covariance estimator formed through
$L_1$-penalized likelihood optimization. We will use a~particular
version to estimate the population weight matrix (equivalently the
inverse error covariance matrix) which was named SPICE by
\citet{Rotetal08} for ``sparse permutation invariant covariance
estimator.''

Let $\widetilde n = n-r-1$ and let
%
%
\begin{equation}\label{Deltahat}
\Deltabfhat= \widetilde{n}^{-1}\Zbb^{T}\Qbf_{\Fbb}\Zbb\in\real
{p \times p}
\end{equation}
denote the\vspace*{1pt} residual sample covariance matrix from the linear regression
of~$\X$ on $\f$.
We construct $\What_{\lambda} = \Deltabfhat{}^{-1}_{\lambda}$ using
%
%
\begin{eqnarray} \label{optomega2}
\Omegabfhat_{\lambda} &\,{=}\,& \mathop{\arg\min}_{\Omegabfs> 0} \biggl[ \tr\{
\diag^{-1/2}(\Deltabfhat)
\Deltabfhat\diag^{-1/2}(\Deltabfhat)
\Omegabf\}\!-\!{\log}|\Omegabf|\!+\!\lambda\sum_{i\neq j} |\Omega_{ij}|\biggr],
\nonumber\hspace*{-40pt}
\\[-8pt]\\[-8pt]
\Deltabfhat{}^{-1}_{\lambda} &\,{=}\,& \diag^{-1/2}(\Deltabfhat)\Omegabfhat_{\lambda}
\diag^{-1/2}(\Deltabfhat),\nonumber\hspace*{-40pt}
\end{eqnarray}
where $\lambda\geq0$ is a~tuning parameter and $\Omega_{ij}$ is the
$(i,j)$ element of $\Omegabf$. The optimization in~(\ref{optomega2})
produces a~sparse estimator $\Omegabfhat_{\lambda}$ of the
inverse error correlation matrix $\Omegabf$, which is then rescaled by the
residual sample standard deviations to
give the inverse error covariance estimator $\Deltabfhat{}^{-1}_{\lambda}$.
The use of $L_1$-penalized likelihood optimization to estimate a~sparse inverse
covariance matrix has been studied
extensively in the literature [\citet{dAsBanElG08},
\citet{YuaLin07}, \citet{FriHasTib08},
\citet{Rotetal08}, \citet{LamFan09}]. Although we impose
sparsity on the off-diagonal entries of
the inverse error covariance matrix $\Deltabf^{-1}$, we expect that in many
situations $\Deltabf^{-1}$ may not be sparse.
We selected this estimator since it is able to adapt, with its tuning parameter,
to both sparse and less sparse estimates, and can lead to a~substantial
reduction in variability of our reduction estimator when $p$ is large; however,
using small values of the tuning parameter $\lambda$ to give
less-sparse inverse
error covariance estimates leads to slow convergence of algorithms to solve
(\ref{optomega2}) when $p$ is large.

\citet{Rotetal08} established the consistency of
$\Deltabfhat{}^{-1}_{\lambda}$ in a~special case of model~(\ref{vecmodel})
characterized by the following four conditions: (A) $r=0$,
so $\Deltabfhat= \Zbb^{T}\Zbb/\tn$ is the usual estimator of the marginal
covariance matrix of~$\X$. (B) The errors $\varepsilonbf$ are normally
distributed with mean 0 and variance~$\Deltabf$. (C)~The largest
$\varphi_{\max}(\Deltabf)$ and smallest $\varphi_{\min}(\Deltabf
)$ eigenvalues
of $\Deltabf$ are uniformly bounded; that is,
for all $p$
%
%
\begin{equation} \label{eqeigen}
0 < \underline{k} \leq\varphi_{\min} (\Deltabf) \leq\varphi
_{\max}
(\Deltabf) \leq\overline{k} < \infty,
\end{equation}
where $\underline{k}$ and $\overline{k}$ are constants. (D) $\lambda
\asymp
\sqrt{\frac{\log p}{ n}}$.
Then [\citet{Rotetal08}]
%
%
\begin{equation} \label{boundspice}
\|\Deltabfhat{}^{-1}_{\lambda} - \Deltabf^{-1}\| = O_p \bigl(\sqrt
{n^{-1}\bigl(s(p)+1\bigr)
\log p}\bigr) ,
\end{equation}
where $s(p)$ is the total number of nonzero off-diagonal entries of
$\Deltabf^{-1}$, which could grow with $p$.

In this application $\Sbf= \W^{-1/2}(\What-\W)\W^{-1/2} =
\Deltabf^{1/2} (\Deltabfhat{}^{-1}_{\lambda} - \Deltabf^{-1} )
\Deltabf^{1/2}$. To find the order of $\|\Sbf\|$ and thereby allow
application of Proposition~\ref{propWhat1} with SPICE, we relax the
conditions of Rothman et al. (\citeyear{Rotetal08}) by allowing (A$^*$) $r
> 0$ and assuming that (B$^*$) the elements $\varepsilon_{j}$ of
$\varepsilonbf$ are \textit{uniformly sub-Gaussian random variables}.
That is, we assume there exist positive constants~$K_1$ and~$K_2$ such
that for all $t > 0$ and all $p$
%
%
\begin{equation}
\label{subgaussian}
P(|\varepsilon_{j} | > t ) \leq K_1 e^{-K_2 t^2},\qquad
j=1,\ldots,p.
\end{equation}
This assumption is compatible with universal conditions (W.1) and (W.2).
When the $\varepsilon_{j}$'s are normal, it is equivalent to requiring that
their variances are uniformly bounded as $p \rightarrow\infty$. The following
lemma gives the orders of~$\|\Sbf\|$ and~$\|\E(\Sbf^{2})\|$ required for
conditions (S.1) and (S.2).
%
%
\begin{lemma}\label{lemmathreshold}
Let $\What= \Deltabfhat{}^{-1}_{\lambda}$, let $s \equiv s(p)$ be the
total number of nonzero off-diagonal entries of $\Deltabf$ and let
$\omegaS(n,p)= \sqrt{n^{-1}(s+1) \log p}$. Under conditions
\textup{(B$^*$)}, \textup{(C)} and \textup{(D)}, $\| \Sbf\| =
O_{p}(\omegaS)$ and $\|\E(\Sbf^{2})\|=O(\omegaS^{2})$.
\end{lemma}

Applying Proposition~\ref{propWhat1}(iii) in the context of SPICE,
$\rhobf= \I_{p}$ because $\W= \Deltabf^{-1}$.
Consequently, $\var(\nubf) = 0$ and
$\psi= \sqrt{p}/h\sqrt{n} \leq\kappa$. From this we immediately
get the
following convergence rate bound for $\Rhat_{\Deltabfshat_\lambda}
\equiv
\Rhatspice$.
%
%
\begin{prop}\label{propWDeltalam}
Assume that the conditions and notation of Lem\-ma~\ref{lemmathreshold} hold.
Then
%
%
\begin{equation} \label{eqspicered}
\Rhatspice(\X_{N}) - \R(\X_{N}) = O_{p} (\kappa) + O_p (\omegaS).
\end{equation}
\end{prop}

If the number $s$ of nonzero off-diagonal entries of $\Deltabf^{-1}$
is bounded
as $ p \rightarrow\infty$ and the regression is abundant ($\kappa
\asymp
n^{-1/2}$),
then $\Rhatspice(\X_{N}) - \R(\X_{N}) = O_{p} (n^{-1/2} \log^{1/2}
p)$. This
allows for the number of variables $p$ to grow much faster than the
sample size,
so long as $s$ is bounded. On the other extreme, if there is no
sparsity so that $s =
p(p-1)$, then the bounding rate implied in Proposition~\ref{propWDeltalam} would
seem to indicate that $n$ needs to be large relative to $\{p(p-1)+1\}
\log p$ for
good results. Based on our results for normal errors in
Section~\ref{secnormalerrors}, we anticipate that this rate is not
sharp, most
notably when~$\Deltabf^{-1}$ is not sparse, and that it can be improved,
particularly when additional structure is imposed (see
Section~\ref{secnormalerrors}).

Conditions (B$^*$), (C) and (D) required for Proposition~\ref{propWDeltalam} are
in addition to the active universal constraints stated in
Section~\ref{secsignal}. Of the universal constraints, we still
require the
signal rate property~(\ref{h}), the order $\kappa= O(1)$
and (W.2). The remaining universal constraints are implied by the
conditions of
the proposition: since $\rhobf= \I_{p}$, the regression is stable by
Lemma~\ref{lemmastable}, (W.1) holds and $\psi\leq\kappa$ so~(\ref{kappa})
holds.

\subsection{A diagonal weight matrix}\label{secdiag}

We include in this section a~discussion of the asymptotic behavior of the
reduction based on the diagonal weight matrix
$\What= \diag^{-1}(\Deltabfhat)$, where $\Deltabfhat$ was defined in
(\ref{Deltahat}). This weight matrix, which corresponds to the population
weight matrix $\W= \diag^{-1}(\Deltabf)$, ignores the residual
correlations and
adjusts only for residual variances. However, in contrast to SPICE,
there is no
tuning parameter involved and so its computation is not an issue.

With $\W= \diag^{-1}(\Deltabf)$, $\rhobf$ is the residual
correlation matrix
and
\[
\Sbf= \diag^{1/2}(\Deltabf)\bigl(\diag^{-1}(\Deltabfhat) - \diag
^{-1}(\Deltabf)\bigr)
\diag^{1/2}(\Deltabf).
\]
The following lemma gives the orders of $\|\Sbf\|$ and $\|\E(\Sbf
^{2})\|$ in
preparation for application of Proposition~\ref{propWhat1}.
%
%
\begin{lemma}\label{lemmadiag}
Let $\What= \diag^{-1}(\Deltabfhat)$, let $\omegaD= n^{-1/2}\log
^{1/2} p$
and assume that the elements $\varepsilon_{j}$ of the errors
$\varepsilonbf$ are
sub-Gaussian random variables. Then $\|\Sbf\| = O_{p}(\omegaD)$ and
$\|\E(\Sbf^{2})\| = O(\omegaD^{2})$.
\end{lemma}

In contrast to Lemma~\ref{lemmathreshold}, here we require only that
the errors
are sub-Gaussian and not uniformly sub-Gaussian. This is because the diagonal
weight matrix effectively standardizes the variances.
%
%
\begin{prop}\label{propWDeltadiag}
Assume that the conditions and notation of Lem\-ma~\ref{lemmadiag} hold.
Then
%
%
\begin{equation} \label{eqdiagred}
\Rhat_{\mathrm{diag}}(\X_{N}) - \R(\X_{N}) = \nubf+ O_{p} (\kappa) +
O_{p}(\psi)
+ O_p (\omegaD).
\end{equation}
\end{prop}

The order~(\ref{eqdiagred}) for the diagonal weight matrix can be
smaller than,
equal to, or greater than the order~(\ref{eqspicered}) for the SPICE weight
matrix, depending on the underlying structure of the regression.
Using results from the discussion of Section~\ref{secTRWhat}, if $\|
\rhobf\| =
O(\sqrt{h})$ and $\spn(\bolds{\gamma})$ is a~reducing subspace
of $\rhobf$, then
$\Rhat_{\mathrm{diag}}(\X_{N}) - \R(\X_{N}) = O_{p}(\kappa) + O_p
(\omegaD)$,
and thus the diagonal weight matrix can produce a~better rate because
$\omegaD
\leq\omegaS$.
If $\Deltabf$ is itself a~diagonal matrix, then $s=0$, $\var(\nubf)
= 0$, $\kappa
\geq\psi$
and $\omegaD= \omegaS$, and consequently
%
%
\begin{equation}\label{diagred2}
\Rhat_{\mathrm{diag}}(\X_{N}) - \R(\X_{N}) = O_{p} (\kappa) + O_p
(\omegaD).
\end{equation}
In this case the two weight matrices result in the same order, so there
seems to be no asymptotic loss incurred by SPICE when $\Deltabf$ is
diagonal.

\section{\texorpdfstring{Normal errors and $\xi= \beta\f$}{Normal errors and xi = beta f}}
\label{secnormalerrors}

We consider in this section the relatively ideal situation in which the
errors $\varepsilonbf$ are normally distributed
and the user-selected basis function $\f$ correctly models the true
coordinates, so $\bolds{\xi}_{i} = \betabf\f_{i}$, $i=1,\ldots
,n$, for some fixed
matrix $\betabf\in\real{d \times r}$ of rank $d$. The diagonal
weight matrix is revisited under these assumptions in
Section~\ref{normalerrorsdiag}. In Section~\ref{normalerrorsn} we add the
sample size constraint $n > p+r+4$ so $\Deltabfhat> 0$. The importance
of this setting is that
we can obtain the oracle rate.

\subsection{Diagonal weight matrix}\label{normalerrorsdiag}

The rates~(\ref{diagred2}) for a~diagonal $\Deltabf$ can be sharpened
considerably when the errors are normal and $\bolds{\xi}= \betabf
\f$:
%
%
\begin{prop}\label{propnormdiag}
Assume that $\varepsilonbf\sim N(0, \Deltabf)$, $\Deltabf$ is a~diagonal matrix, $\bolds{\xi}= \betabf\f$, $n > r+5$ and $\What=
\diag^{-1}(\Deltabfhat)$.
Then
$\Rd(\X_{N}) - \R(\X_{N}) = O_{p} (\kappa)$.
\end{prop}

We see from this proposition that when the errors are normal, the order
$O_{p}(\omegaD)$ that appears in~(\ref{diagred2}) is no longer present
and the rate for $\Rd$ reduces to the oracle rate $\kappa^{-1}$. The
condition $n > r+ 5$ stated in Proposition~\ref{propnormdiag} is
needed to insure
the existence of the variance of an inverse chi-squared random variable.
Its proof is omitted from the supplement since it follows the same
lines as
the proof of Proposition~\ref{propWhat1}.

\subsection{$n > p + r + 4$}\label{normalerrorsn}

Recalling that $\Deltabfhat\in\real{p \times p}$ is the residual
covariance matrix defined in~(\ref{Deltahat}), this constraint on $n$
means that $\Deltabfhat> 0$ with probability 1, allowing the
straightforward use of $\What= \Deltabfhat ^{-1}$ as the weight
matrix.\vspace*{1pt} The population weight matrix is $\W= \Deltabf
^{-1}$ with corresponding $\Sbf= \Deltabf^{1/2} \Deltabfhat{}^{-1}
\Deltabf^{1/2} - \I_{p}$. The asymptotic behavior of $\|\Sbf\|$ in this
setting can be obtained as follows: phrased in the context of this
article, suppose that $r = 0$ so $\Deltabfhat=
(n-1)^{-1}\sum_{i=1}^{n}(\X_{i}-\Xbar)(\X_{i}- \Xbar)^{T}$. Johnstone
(\citeyear{Joh01}) [see also Bai (\citeyear{Bai99}), page
635]\vspace*{1pt} showed that\vadjust{\goodbreak} $\varphi_{\max}( \Deltabf^{-1/2}
\Deltabfhat\Deltabf^{-1/2}) - 1 \asymp_p (\sqrt{p/n})$ and
\citet{Pau05} showed that\vspace*{1pt} $\varphi_{\min} ( \Deltabf^{-1/2}
\Deltabfhat\Deltabf^{-1/2}) - 1\asymp_p (\sqrt{p/n})$ when $p/n
\rightarrow a~\in[0,1)$. Together these results imply that $\| \Sbf\|
\asymp_p (\sqrt{p/n})$ and therefore $\| \Sbf \| \rightarrow0$ in
probability if and only if $p/n \rightarrow0$, which gives the order
for condition (S.1). Still with $r=0$, since $(n-1) \Deltabfhat{}^{-1}$
is distributed as the inverse of a~Wishart $W_{p}(\Deltabf, n-1)$
matrix, we used results from \citet{von88} to verify condition
(S.2). Combining these results with Proposition~\ref{propWhat1} we have
%
%
\begin{equation}\label{normrate1}
\Rg(\X_{N}) - \R(\X_{N}) = O_{p}\bigl(\sqrt{p/n}\bigr),
\end{equation}
where $\Rg$ denotes the reduction with $\What= \Deltabfhat{}^{-1}$.
This suggests that it may be reasonable to use $\Deltabfhat{}^{-1}$ as the
weight matrix only when $n$ is large relative to $p$, essentially
sending us
back to the usual requirement. In the next proposition we show this is not
the case with $r > 0$ because the order in~(\ref{normrate1}) is in
fact too
large.
%
%
\begin{prop}\label{propWDelta1}
Assume\vspace*{1pt} that $\varepsilonbf\sim N_{p}(0,\Deltabf)$, $\bolds{\xi}=
\betabf\f$, $n > p + r + 4$ and
$p/n \rightarrow[0,1)$. Let $a~= (n-p-1)^{-1}$. Then when $\What=
\Deltabfhat{}^{-1}$:
\begin{longlist}
\item
$
h^{-1}\E(\Khat) = a(n-r-1)\{\K+ \kappa^2
\I_r \},
$
\item $h^{-1}(p)\{\Khat- \E(\Khat)\} = O_{p}(1/\sqrt{n})$,
\item $\Rg(\X_{N}) - \R(\X_{N}) = O_{p}(\kappa)$.
\end{longlist}
\end{prop}

Recall from Corollary~\ref{corW21} that when $\What= \Deltabf ^{-1}$,
we obtain the oracle rate $\Rhat_{\Deltabfs}(\X_{N}) - \R(\X_{N}) =
O_{p}(\kappa)$. Comparing this result with the conclusion of
Proposition~\ref{propWDelta1}(iii), we see that there is no impact on
the rate of convergence when using $\What= \Deltabfhat{}^{-1}$ instead
of $\What= \Deltabf^{-1}$, both weight matrices giving estimators that
converge at the oracle rate. In contrast to the gross bound given in
(\ref{normrate1}), Proposition~\ref{propWDelta1} indicates that it may
in fact be quite reasonable to use $\Deltabfhat{}^{-1}$ as the weight
matrix when $n > p + r + 4$, without requiring $n$ to be large relative
to $p$.

Properties of the normal and the accurate modeling of $\bolds{\xi
}= \betabf\f$
surely contribute to the oracle rate of
Proposition~\ref{propWDelta1}. Recall from~(\ref{Khat}) that $\Khat=
\Fd^{1/2}\Bhat^{T}\Whati\Bhat\Fd^{1/2}$. Under the conditions of
Proposition~\ref{propWDelta1}, $\Bhat\indep\Whati$, which is not
required for
our most general results given in Proposition~\ref{propWhat1}. Additionally,
we were able to use exact calculations in places
where bounding was otherwise employed. The inverse residual sum of squares
matrix $(n-r-1) \Deltabfhat{}^{-1}$ is distributed as the inverse of a~Wishart
$W_{p}(\Deltabf, n-r-1)$ matrix and the moments of an inverse Wishart
derived by
\citet{von88} were used extensively in the proof of
Proposition~\ref{propWDelta1}. For instance, although~$\Deltabfhat$
has an
inverse with probability 1 when $n > p + r + 1$, the constraint on $n$
in the
hypothesis is needed to ensure the existence of the variance of the inverse
Wishart.\vadjust{\goodbreak}

The signal rate~(\ref{h}) and the order $\kappa= O(1)$ are still
required for Proposition~\ref{propWDelta1}, but its hypothesis implies
all the other universal conditions stated in Section~\ref{secsignal}:
since $\W= \Deltabf^{-1}$ we have $\rhobf= \I_{p}$ and thus (1) the
regression is stable by Lemma~\ref{lemmastable}, (2) $\psi\leq\kappa$
so~(\ref{kappa}) holds, (3) $\E(\varepsilonbf^{T}\W\varepsilonbf) = p$
so condition (W.1) holds, and (4) condition (W.2) holds by Lemma
\ref{lemmaW}(i).

\section{Simulations}\label{secsimulations}

The effects of predictor correlations are often a~concern when dealing
with high-dimensional regressions. In this section, we introduce
predictor correlations into our simulations to give some intuition
about their impact on the four reduction estimators.

The computation\vspace*{1pt} of the four reduction estimators $\Ri,
\Rhatspice, \Rd$ and~$\Rg$ relies on the computation of their weight
matrix estimators $\What$, which are available in closed form for $\Ri,
\Rd$ and $\Rg$; however, computing the weight matrix estimator for
$\Rhatspice$ requires us to select its tuning parameter~$\lambda$ and
to use an iterative algorithm.

\subsection{\texorpdfstring{Computation and tuning parameter selection for $\Rhatspice$}{Computation and tuning parameter selection for R spice}}
Several computational algorithms [\citet{dAsBanElG08},
\citet{YuaLin07}, \citet{FriHasTib08}, \citet{Rotetal08}]
have been proposed to compute the solution to~(\ref{optomega2}),
of which we propose to use the graphical lasso (glasso) algorithm of
\citet{FriHasTib08}.
We employ $K$-fold cross-validation to select the tuning parameter
$\lambda$ for
the weight matrix estimator $\What_{\lambda} \equiv
\Deltabfhat{}^{-1}_{\lambda}$,
where we minimize the validation negative log-likelihood.
Specifically, we solve $\hat\lambda= \arg\min_{\lambda} \sum_{k=1}^{K}
g_k(\lambda)$
where
\begin{eqnarray*}
g_k(\lambda) & = & \tr\bigl[n_{k}^{-1}\A_{(k, \lambda)}^{T}\A_{(k,
\lambda)}
\What_{\lambda}^{(-k)}
\bigr] - {\log}\bigl|\What_{\lambda}^{(-k)}\bigr|, \\
\A_{(k, \lambda)} & = & {\Zbb}^{(k)} - {\Fbb}^{(k)} \bbfhat^{(-k)
T}_{\lambda} \Gammabfhat{}^{(-k) T}_{\lambda}.
\end{eqnarray*}
In the above expression, a~superscript of $(k)$ indicates that the
quantity is based on the observations inside the $k$th fold, and
a~superscript of $(-k)$ indicates that the quantity is based on observations
outside of the $k$th fold, and $n_k$ is
the number of observations in the $k$th fold.

\subsection{Simulation settings}
The simulated data were generated from $ \X_i = \Gammabf\bolds{\xi}_i +
\varepsilonbf_i$, $i=1,\ldots,n$, where $Y_1, \ldots, Y_n$ is a~sequence of independent random variables, $\varepsilonbf_1, \ldots,
\varepsilonbf_n$ are independent copies of $\varepsilonbf$, a~$p$-variate normal random vector with mean 0 and variance $\Deltabf$,
and $\bolds{\xi}_{i} = \bolds{\xi}(Y_{i})$ is specified later. As in
the second simulation of Section~\ref{secsimW}, $\Gammabf$ and
$\Deltabf$ were produced anew at the start of each replication, with
the elements of $\Gammabf\in\real{p}$ always generated as a~sequence of
independent standard normal variables and $\Deltabf$ generated as
$\D^{1/2}\Thetabf\D^{1/2}$ where $\D$ is a~diagonal matrix with
diagonal elements sampled from a~uniform\vadjust{\goodbreak} $(1,101)$ distribution and
$\Thetabf= (\theta_{ij})$ is a~correlation matrix with exponentially
decreasing correlations, $\theta_{ij} = \theta^{|i-j|}$ with parameter
$0 < \theta< 1$. This structure has been used frequently to assess the
impact of predictor correlations on high-dimensional methodology [see,
e.g., \citet{BicLev08N2}, \citet{LiYin08}]. For all four
reduction estimators, $h\asymp p$ implying $\kappa= n^{-1/2}$.

For each replication, we determined $\Rhat(\X_{N,j})$ and $\R(\X_{N,j})$
at $j=1,\ldots,100$ new data points generated from the original model. We
assess performance by computing the
magnitude of the sample correlation between $\Rhat$ and $\R$ over the 100
new
data points.
The results are based on 200 independent replications.

For a~given weight matrix $\W$, the expected
signal strength over simulation replications
is determined by $\Es(\Gammabf^T \W\Gammabf) = \Es\tr(\W)$,
where $\Es$
denotes expectation over the simulations.
For $\W= \diag^{-1}(\Deltabf)$, $\Es[\tr(\W)] = \log(101)/100$,
but for
$\W= \Deltabf^{-1}$,
\[
\Es[\tr(\W)] = \tr(\Thetabf^{-1}) \log(101)/100= \frac
{2+(p-2)(1+\theta^2)}{1-\theta^2} \frac{\log
101}{100}.
\]
From this we see that the expected signal strength increases with
$\theta$ when using
$\Rg$ and $\Rhatspice$, but is constant in $\theta$ when using $\Rd
$ and $\Ri$.
In addition, $\|\Deltabf\| = O(1)$ since $\|\Deltabf\| \leq
\|\Deltabf\|_{\infty} = \max_{j}\sum_{i}\delta_{ij} \leq2\cdot
101/(1-\theta)$,
where $\delta_{ij}$ is element $(i,j)$ of $\Deltabf$.

With this general setup we next report results for various combinations of~$n$,
$p$, $\bolds{\xi}$ and $\f$. We extended the applicability
of $\Rg$ to
regressions in which $n < p$ by using
the Moore--Penrose generalized inverse of $\Deltabfhat$, although we have
presented no asymptotic results for this case.

\subsection{\texorpdfstring{Correctly specified $\xi= \beta\f$}{Correctly specified xi = beta f}}

In this section we consider a~simple case where $\bolds{\xi}= Y$,
giving $d=1$, and
$\f= (Y, Y^2, Y^3, Y^4)$, so $r = 4$ and $\bolds{\xi}= \betabf
\f$, where
$\betabf= (1,0,0,0)$ and $Y$ was generated as a~standard normal variate.

\subsubsection{\texorpdfstring{$n=p/2$ and $p \rightarrow\infty$}{n = p/2 and p -> infinity}}\label{simnp1}

In this setting, $n$ and $p$ grow with $n=p/2$.
For the reduction estimators $\Ri$ and $\Rd$, $\|\rhobf\| = O(1)$,
$\|\rhobf\|_F \leq\sqrt{p} \|\rhobf\| = O(\sqrt{p})$ which imply
that $\psi=
\kappa= n^{-1/2}$.
Also, $\var(\nubf) \leq\G_{h}^{-1}\|\Deltabf\|/h = O(p^{-1})$,
hence $\nubf= O_{p}(p^{-1})$. Thus in this setting with
$n \asymp p$, $\Ri$ is $\sqrt{n}$-consistent and~$\Rd$ is at least
$\sqrt{n/\log n}$-consistent as $n,p \rightarrow
\infty$.

Although $\|\Deltabf\| = O(1)$ and $\Deltabf^{-1}$ is a~tri-diagonal matrix,
our theoretical bounds for $\Rhatspice$ guarantee consistency for this model
only when $p$ is bounded and
$n \rightarrow\infty$; however, a~result established by \citet
{Ravetal11}, which requires additional
assumptions, indicates that the weight matrix estimator for $\Rhatspice
$ is
consistent when $\Deltabf^{-1}$ is
tri-diagonal.

%
\begin{figure}
\begin{tabular}{@{}cc@{}}

\includegraphics{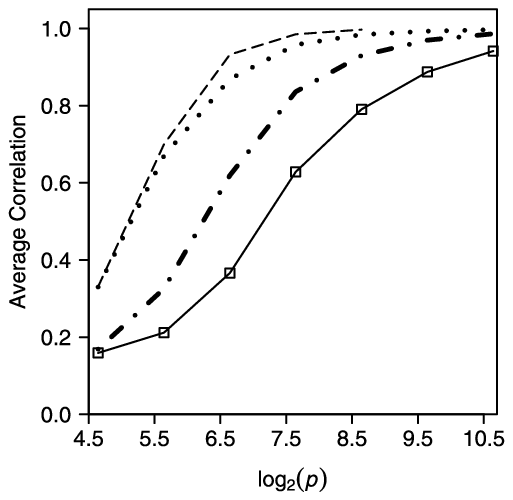}
 & \includegraphics{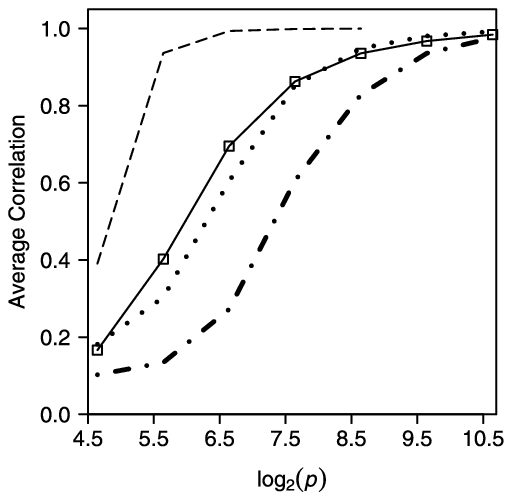}\\
(a) & (b)\\[6pt]
\multicolumn{2}{@{}c@{}}{
\includegraphics{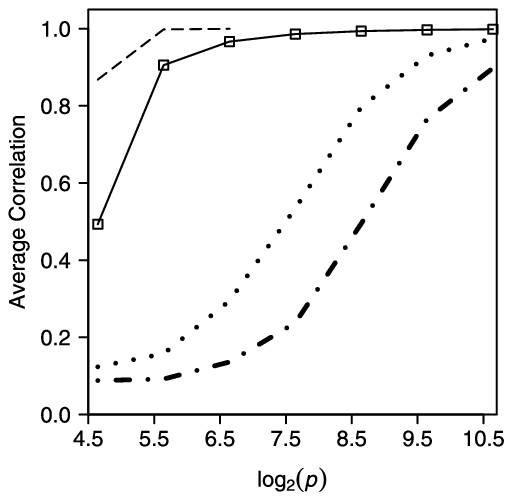}
}\\
\multicolumn{2}{@{}c@{}}{(c)}\vspace*{-3pt}
\end{tabular}
\caption{Comparison of the four estimators of $\R$: $\Rhatspice$
(dashes), $\Rg$ (solid), $\Rd$ (dots) and $\Ri$ (dash dot), with
exponential error correlations and $n=p/2$. Here
\textup{(a)}~$\theta_{ij} =
0.5^{|i-j|}$; \textup{(b)}~$\theta_{ij} = 0.9^{|i-j|}$; \textup{(c)}
$\theta_{ij} = 0.99^{|i-j|}$.}\label{figpexp}
\end{figure}

The results for $p$ and $n$ growing with $n=p/2$ are shown in
Figure~\ref{figpexp}(a)--(c). All reduction estimators
appear to be converging to the population reduction as\vadjust{\goodbreak} $n, p
\rightarrow
\infty$,
even though consistency is not guaranteed in this setting for
$\Rhatspice$ and
$\Rg$,
%
%
\begin{figure}
\begin{tabular}{@{}cc@{}}

\includegraphics{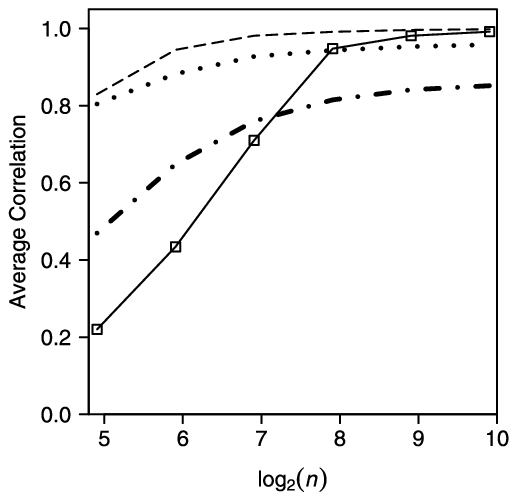}
 & \includegraphics{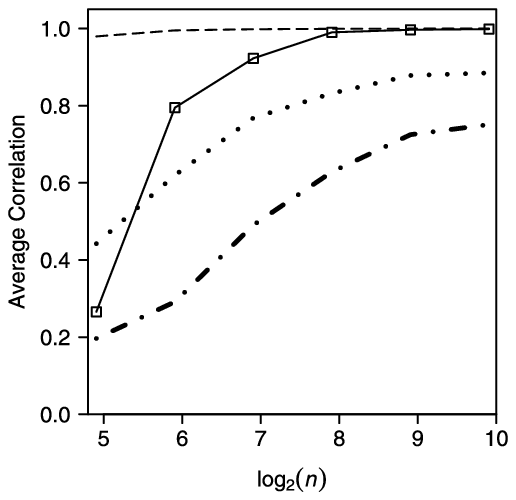}\\
(a) & (b)\\[4pt]
\multicolumn{2}{@{}c@{}}{
\includegraphics{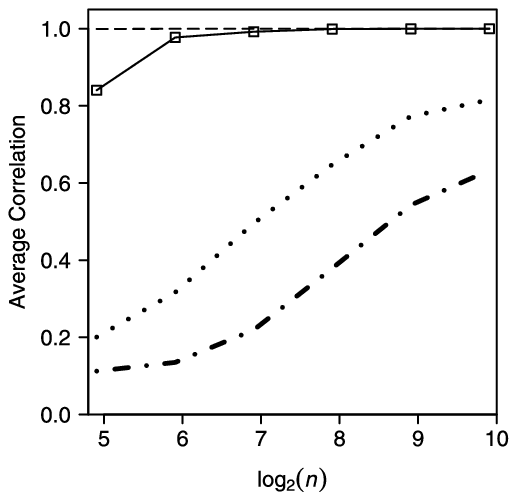}
}\\
\multicolumn{2}{@{}c@{}}{(c)}\vspace*{-3pt}
\end{tabular}
\caption{Comparison of the four estimators of $\R$:
$\Rhatspice$ (dashes), $\Rg$ (solid), $\Rd$ (dots)
and $\Ri$ (dash dot), with exponential error correlations
and $p=100$. Here \textup{(a)}~$\theta_{ij} = 0.5^{|i-j|}$;
\textup{(b)}~$\theta_{ij} = 0.9^{|i-j|}$;
\textup{(c)} $\theta_{ij} = 0.99^{|i-j|}$.}\label{fignexp}\vspace*{-3pt}
\end{figure}
which uses a~Moore--Penrose generalized inverse of the residual sample
covariance matrix. Interestingly, when $p > n$, $\Rg$ outperforms $\Rd$
and $\Ri$ when $\theta\geq 0.9$. Results for $\Rhatspice$ were computed
up to $p=400$ when $\theta\leq0.9$ and up to $p=100$ when $\theta=0.99$
due to intractable computation time required for the glasso algorithm.
In scenarios when $\Rhatspice$ was computed, it outperformed the other
reduction estimators, particularly when $\theta=0.9$. As expected,
larger values of $\theta$ lead to favorable performance for the
reduction estimators with population weight matrix $\W= \Deltabf^{-1}$.

\subsubsection{\texorpdfstring{$p=100$ and $n \rightarrow\infty$}{p = 100 and n -> infinity}}\label{simnp2}
In this setting we fix $p=100$ and let $n$ grow.
Our theory guarantees that $\Rhatspice$ and $\Rg$ are both
$\sqrt{n}$-consistent.
On the other hand, $\Ri$ and $\Rd$ are inconsistent
since $\spn(\bolds{\gamma})$ is not a~reducing subspace of
$\rhobf$ and $p$ is
bounded,
implying $\nubf$ fails to vanish.\vadjust{\goodbreak}

The results for $p=100$ and $n$ growing are illustrated in
\mbox{Figure~\ref{fignexp}(a)--(c)}.
As our theory suggests,
the reduction estimators $\Rhatspice$ and $\Rg$ appear to be converging
to the population reduction as $n$ increases. We see that~$\Rhatspice$
outperformed the other reduction estimators, particularly for
relatively\break
small~$n$ when $\theta=0.9$.
As expected, $\Rg$ outperforms $\Rd$ and~$\Ri$ when~$n$ is much
larger than $p$
or when $\theta$ is large.\vspace*{-3pt}

\subsection{\texorpdfstring{Results for $\xi\neq\beta\f$}{Results for xi not= beta f}} \label{simnp3}

\begin{figure}
\begin{tabular}{@{}cc@{}}

\includegraphics{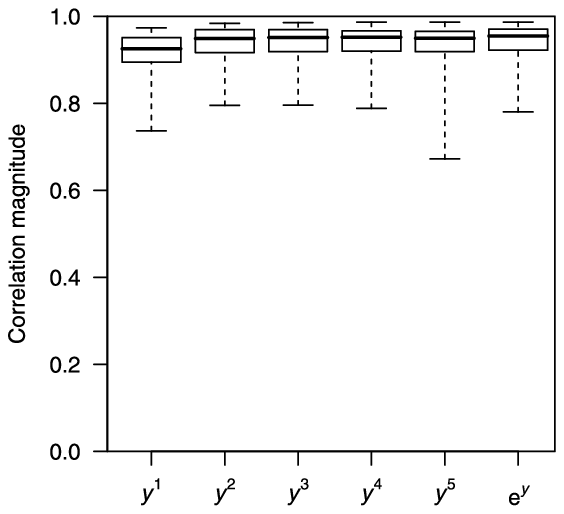}
 & \includegraphics{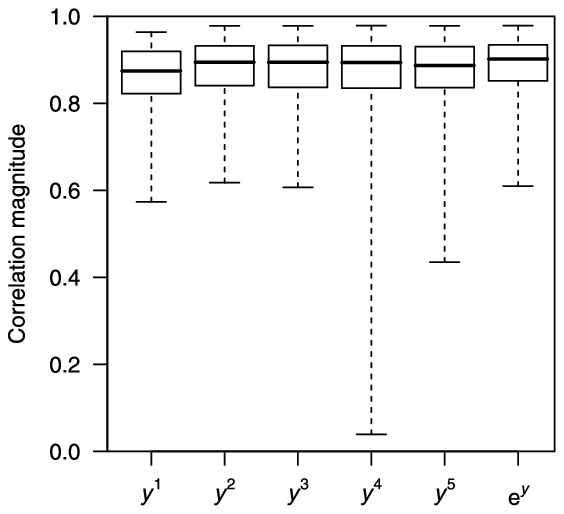}\\
(a) & (b)\\[6pt]

\includegraphics{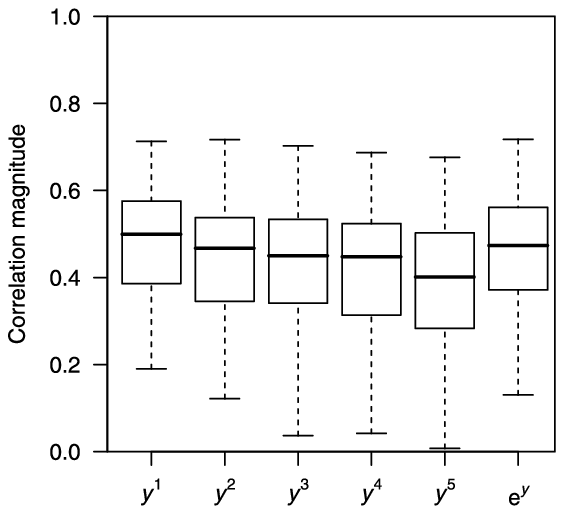}
 & \includegraphics{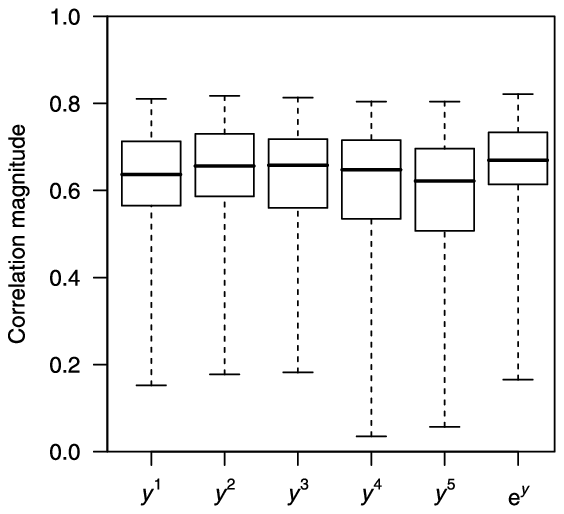}\\
(c) & (d)
\end{tabular}
\caption{Estimators of $\R$ when $\f$ is misspecified, using
exponential error correlations: $\theta_{ij} = 0.5^{|i-j|}$, $n=50$
and $p=100$. Boxplots are labeled by the highest order term in $\f$.
Here \textup{(a)} $\Rhatspice$; \textup{(b)} $\Rd$; \textup{(c)} $\Rg$;
\textup{(d)} $\Ri$.}\label{figbasis}
\end{figure}

In this section we present results for a~misspecified $\bolds{\xi }$
using $\bolds{\xi}= \var^{-1/2}(\exp(Y))[\exp(Y) - \E(\exp(Y))]$ where
$Y \sim\operatorname{Unif}(0,4)$. Holding $n=50$ and $p=100$, we varied
$\f= (y, y^2, \ldots, y^k)^T$ for $k=1,\ldots, 5$ and $\f= \exp(y)$.
The results are summarized in Figure~\ref{figbasis} as
%
%
boxplots of the correlation magnitudes over the 200 replications. Most
striking is how little the choice of $\f$ seems to affect the
estimators. This is likely because there are fairly strong\vadjust{\goodbreak} correlations
between $\bolds{\xi}$ and its approximations provided by the six $\f$'s
used in the simulation, which satisfies condition~(\ref{rank}). The
relationship between the estimators are similar to those in Figure
\ref{fignexp} at $n=50$, regardless of the choice of $\f$.

\section{Spectroscopy application}\label{secdata}

$\!\!$To illustrate operation of the proposed meth\-odology, we examine data from
the potential application area mentioned in the
\hyperref[secintro]{Introduction}. The response
variable is the percentage of fat in samples of beef or pork, and the predictors
are the absorbance spectra ($\log(1/R)$) from near-infrared
transmittance for
fat measured at every second wavelength between 850 and 1050 nm, giving $p=100$
with $n=54$ [\citet{Sbetal07}]. The goal of the study is to predict
the response from the absorbance spectra of a~new sample. Our predictive
framework is described in the next section.
We return to the spectroscopy application in Section~\ref{secspectroscopy}.\vadjust{\goodbreak}

\subsection{Prediction}\label{secprediction}

We predict an unobserved response $Y_{N}$ associated with a~new
observed vector
of predictors $\X_{N}$ by using the forward regression mean function:
$Y_{\pred}
= \E\{Y|\X_{N}\} = \E\{Y|\R(\X_{N})\}$. However, the reduction~$\R
(\X)$ was
based on the inverse regression $\X|Y$ and the development did not
produce a~direct estimator of $\E\{Y|\R(\X_{N})\}$. There are perhaps several
ways of
using an estimated reduction for predicting a~new response. Some
authors have
used standard data-analytic methods to develop predictive models based
on $Y$
and the estimated reduction, and there are a~variety of nonparametric
regression methods that could also be used as well. We follow
\citet{AdrCoo09} and use a~kernel-type estimator of $\E\{Y|\R(\X
_{N})\}$ based
on the
relationship
$ \E\{Y|\X\,{=}\,\bx\}
\,{=}\, \E\{Y|\R(\bx)\}
\,{=}\, \E\{Y g(\R(\bx)|Y)\}/\E\{g(\R(\bx)|Y)\},
$
where $g$ is the conditional density of $\R|Y$. This provides a~method to
estimate $\E\{Y|\X\}$:
%
%
\begin{eqnarray}\label{Ehat}
\widehat{\E}\{Y|\X=\bx\} & = &
\sum_{i=1}^n w_i(\bx) Y_i, \nonumber\\[-8pt]\\[-8pt]
w_i(\bx) & = & \frac{\ghat(\Rhat(\bx)|Y_i)}{\sum_{i=1}^n
\ghat(\Rhat(\bx)|Y_i)}, \nonumber
\end{eqnarray}
where $\ghat$ denotes an estimate of the density and $\Rhat$ is
the estimated reduction. This estimator is reminiscent of a~nonparametric kernel estimator, but there are consequential
differences. The weights in a~kernel estimator do not depend on the
response, while the weights $w_i$ here do. Kernel weights typically
depend on the full vector of predictors $\X$, while the weights here
depend on $\X$ only through the~estimated reduction $\Rhat(\bx)$.
Multivariate kernels are usually taken to be the product of univariate
kernels, corresponding here to treating the components of $\R$ as
independent. Finally, there is no need for bandwidth estimation because
the weights are determined entirely from $\ghat$.

When $d$ is small relative to $p$ it may often be reasonable to assume that
$\R(\X)|Y$ is normally distributed, which seems appropriate for the
spectroscopy
data. Ignoring constants not depending on $i$, we have
\[
g(\R({\mathbf x})|Y_{i}) \propto\exp\bigl\{-(1/2)\bigl(\R({\mathbf x})
- \bolds{\xi}_{i}
\bigr)^{T}(\Gammabf^{T}\Deltabf^{-1}\Gammabf)\bigl(\R({\mathbf x}) -
\bolds{\xi}_{i}\bigr)\bigr\}.
\]
Substituting the estimators $\Gammabfhat$, $\bbfhat\f_{i}$ and
$\What$ for
$\Gammabf$, $\bolds{\xi}_{i}$ and $\Deltabf^{-1}$ gives the
weights required for
(\ref{Ehat}). For $\Ri$, we set $\What= (p/\tr(\Deltabfhat)) \I_{p}$.

%
\begin{figure}[b]
\begin{tabular}{@{}cc@{}}

\includegraphics{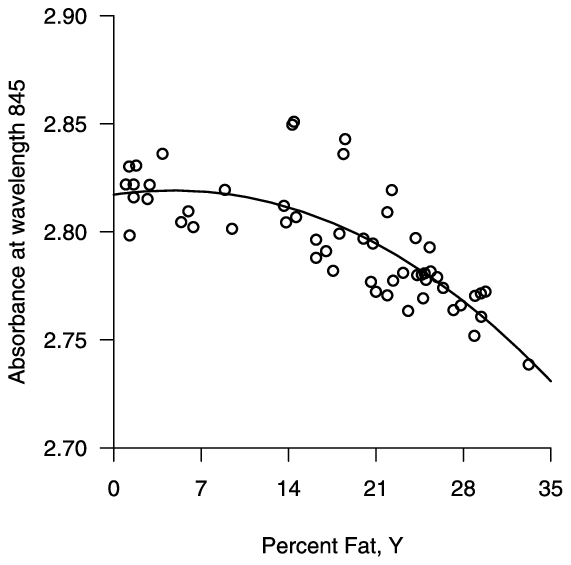}
 & \includegraphics{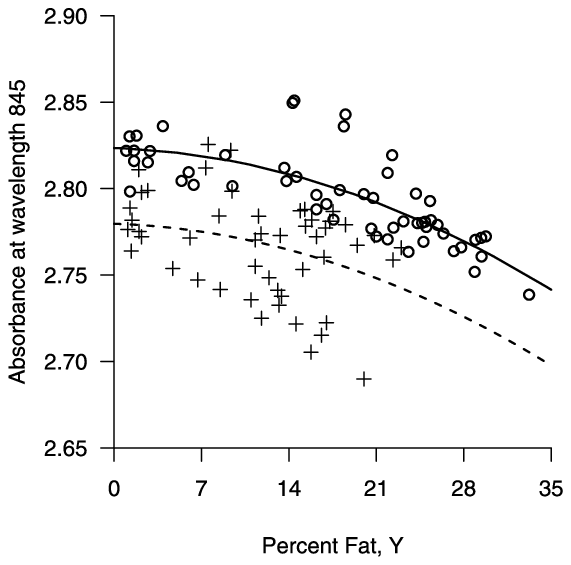}\\
(a) & (b)\vspace*{-3pt}
\end{tabular}
\caption{Inverse response plot of the absorbance at wavelength 854
versus $Y$ for \textup{(a)} pork samples and \textup{(b)} pork and beef
samples. The line for the pork sample is a~fitted cubic polynomial. For
pork and beef the lines represent a~second-order polynomial fit to the
data with one intercept for pork (solid) and a~second intercept for
beef (dashes).} \label{figporkinvresp}\vspace*{-3pt}
\end{figure}

\subsection{Spectroscopy}\label{secspectroscopy}
Since there are only $p=100$ predictors it was straightforward, albeit somewhat
tedious, to inspect inverse response plots of $X_{j}$ versus $Y$,
$j=1,\ldots,100$, to gain information about the likely structure of~$\f$.
We performed two analyses; the first consisted of 54 pork samples and
the second
consisted of 103 meat samples
of beef and pork. For a~specified value of $d$, we assessed\vadjust{\goodbreak} the fitted model
by using the residual mean square $\rms(d,r,\Rhat_{(\cdot)}) =
\sum_{i=1}^{n}(Y_{i} - \Yhat_{i})^{2}/n$, where the fitted values
\mbox{$\Yhat_{i}=\widehat{\E}(Y|\X=\X_{i})$} were determined using~(\ref{Ehat}) with
the indicated combination of $d$, $r$ and reduction~$\Rhat_{(\cdot)}$.

We selected the value of $d$ by adapting the permutation scenario
developed by
\citet{CooYin01}. The hypothesis $d = d_{0}$ was tested sequentially,
starting at $d_{0} =0$ and estimating $d$ as the first hypothesized
value that
was not rejected. The test statistic $\rms(d_{0}+1, r, \R_{(\cdot
)})$ was
compared to the distribution of $\rms(d_{0}+1, r, \R_{(\cdot)})$
induced by
$1000$ random permutations~$\J$ applied to the rows of the predictor matrix
$\Xbb$ as follows:
\[
\Xbb_{\mathrm{perm}} = \Xbb\Pbf_{\Gammabfshat(\What)}^{T} +
\J\Xbb\Qbf_{\Gammabfshat(\What)}^{T},\vspace*{-3pt}
\]
where $\Gammabfhat$ and $\What$ were computed under the null
hypothesis. This
scheme leaves the signal $\Xbb\Pbf_{\Gammabfshat(\What)}^{T}$
intact while
permuting the uninformative part of the predictors $\Xbb\Qbf
_{\Gammabfshat
(\What)}^{T}$.

The multicollinearity of the predictors made computing the weight matrix
estimator for $\Rhatspice$ difficult for
small values of its tuning parameter, values for which our cross-validation
procedure recommended.
We subsequently set its tuning parameter to $\lambda=2^{-10}$ for both
analyses,
since
this was the smallest value for which a~numerically stable solution was
available.

\subsubsection{Analysis of pork samples}
In this case we concluded that a~cubic polynomial $\f(y) = (y, y^{2},
y^{3})^{T}$ would be adequate; a~representative plot is shown in
Figure~\ref{figporkinvresp}(a).
Performing the permutation test to select $d$,
the test statistic for\vadjust{\goodbreak} $d=0$ using $\Rg$ was $\rms(1,3,\Rg) = 0.31$
which was
smaller than the smallest value $10.7$ of $\rms(1,3,\Rg)$ observed
among the
$1000$ random permutations under the hypothesis that $d=0$. Since the test
statistic is much smaller than can be accounted for by chance under the null
hypothesis, we concluded that $d \geq1$. Similarly, to test $d=1$, we observed
$\rms(2,3,\Rg)= 0.29$ which fell at the 80th quantile of the permutation
distribution of $\rms(2,3,\Rg)$ under the hypothesis. Consequently,
we used
$d=1$ for the model. Comparisons with other values of $r$ figured in
our choice
$r=3$. Cross-validation using $\rms$ as the criterion might also be
used to
select $d$.\looseness=-1

%
\begin{figure}[b]
\vspace*{-3pt}
\begin{tabular}{@{}cc@{}}

\includegraphics{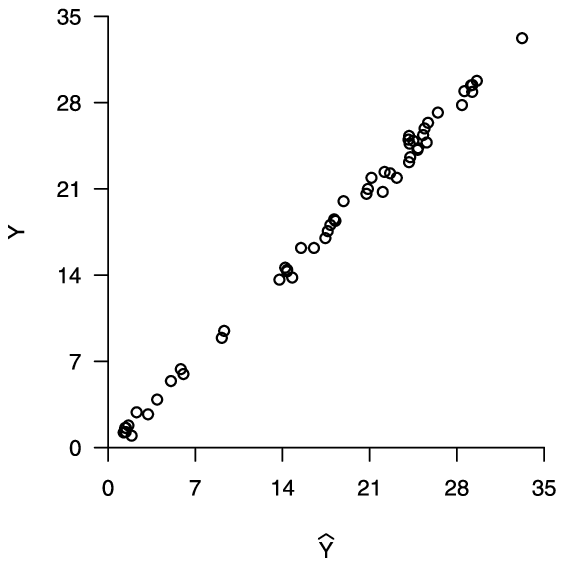}
 & \includegraphics{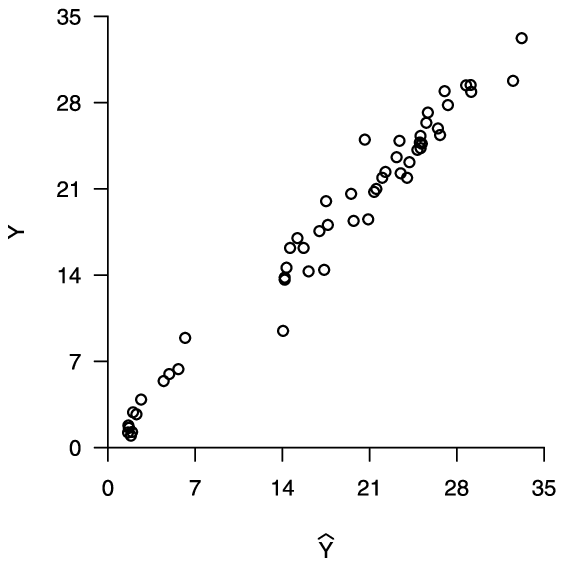}\\
(a) & (b)\\[4pt]

\includegraphics{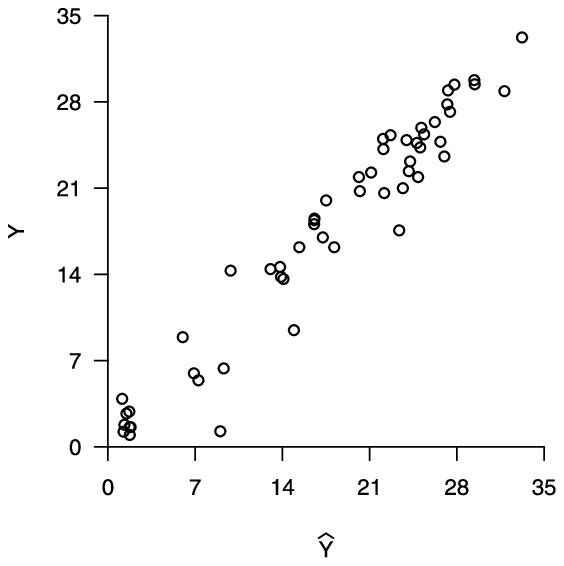}
 & \includegraphics{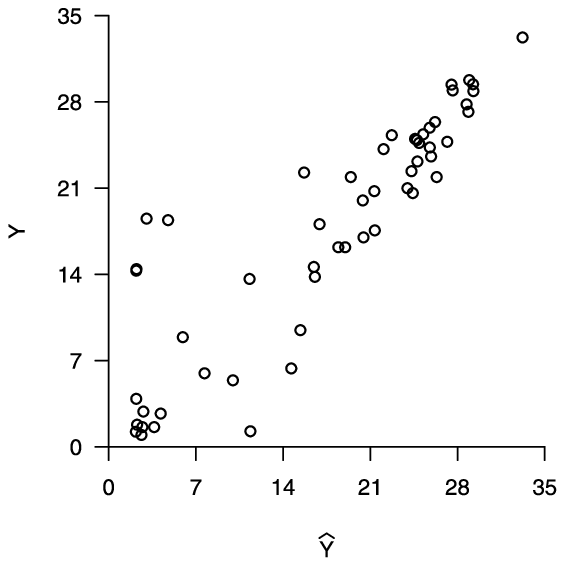}\\
(c) & (d)\vspace*{-3pt}
\end{tabular}
\caption{Plots of the response versus fitted values for the
four estimators \textup{(a)} $\Rg$, \textup{(b)}~$\Rhatspice$,
\textup{(c)}~$\Rd$ and \textup{(d)} $\Ri$,
using the pork samples.}\label{figpork}
\end{figure}

Fitting each of the four estimators with $d = 1$, we found that $\rms
(1,3,\allowbreak\Rhat) = 21.55$, $5.15$, $0.31$ and $2.15$ for $\Rhat= \Ri$,
$\Rd$, $\Rg$ and $\Rhatspice$. Plots of $Y$ versus $\Yhat$ are shown in
Figure~\ref{figpork}(a)--(d). The predictors in this illustration are
highly collinear, as is typical in spectral data, and the relative
performance of the four estimators is qualitatively similar to that
shown in previous\vspace*{1pt} simulations; however, numerical instability degraded
the performance of $\Rhatspice$. The relative\vspace*{1pt} signal rates
for the four estimators are reflected by the values of
$\Gammabfhat{}^{T}\What\Gammabfhat= 35.3$, $128.4$, $169.7$ and $63.38$
for $\What= (p/\tr(\Deltabfhat))\I_{d}$, $\diag^{-1}(\Deltabfhat)$,
$\Deltabfhat{}^{-}$ and $\Deltabfhat{}^{-1}_{\hat\lambda}$.

\subsubsection{Analysis of both pork and beef samples}
In this case we concluded that a~second-order polynomial and the indicator
function of beef
would be adequate, $\f(y) = (y, y^{2}, J(\mathrm{beef}))^{T}$; a~representative
plot is shown in Figure~\ref{figporkinvresp}(b).
Using the permutation test approach to select $d$,
the test statistic for $d=0$ using $\Rg$ was $\rms(1,3,\Rg) = 0.55$
which fell
at the 0.003 quantile of $\rms(1,3,\Rg)$ observed among the $1000$ random
permutations under the hypothesis that $d=0$. Since the test statistic
is much
smaller than can be accounted for by chance under the null hypothesis, we
concluded that $d \geq1$. Similarly, to test $d=1$, we observed
$\rms(2,3,\Rg)= 0.01$ which fell at the 0.70 quantile of the permutation
distribution of $\rms(2,3,\Rg)$ under the hypothesis. Consequently,
we used
$d=1$ for the model.

%
\begin{figure}
\begin{tabular}{@{}cc@{}}

\includegraphics{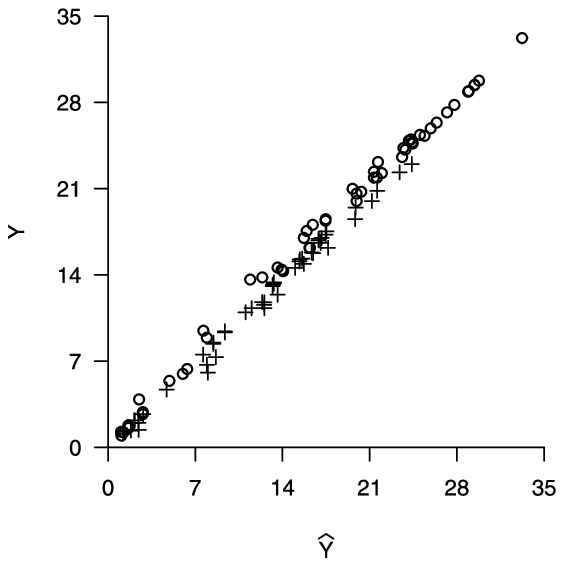}
 & \includegraphics{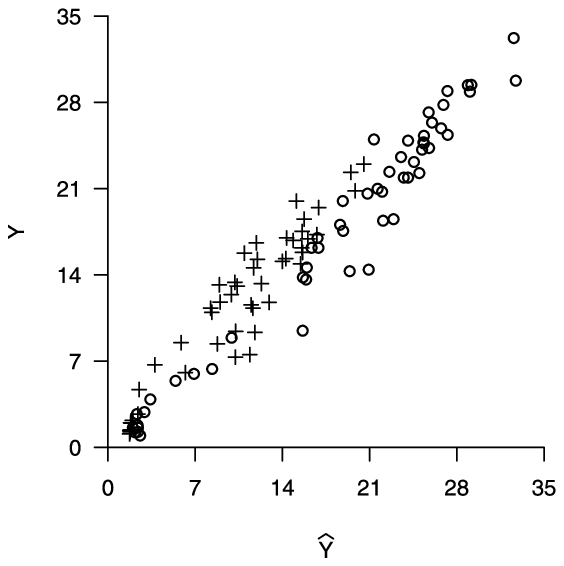}\\
(a) & (b)\\[6pt]

\includegraphics{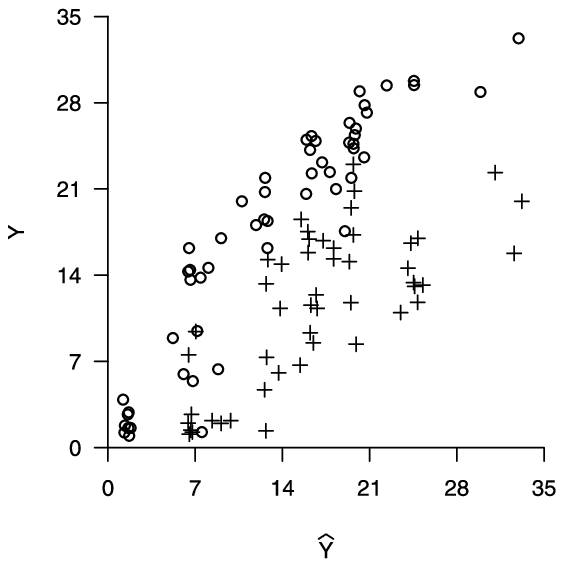}
 & \includegraphics{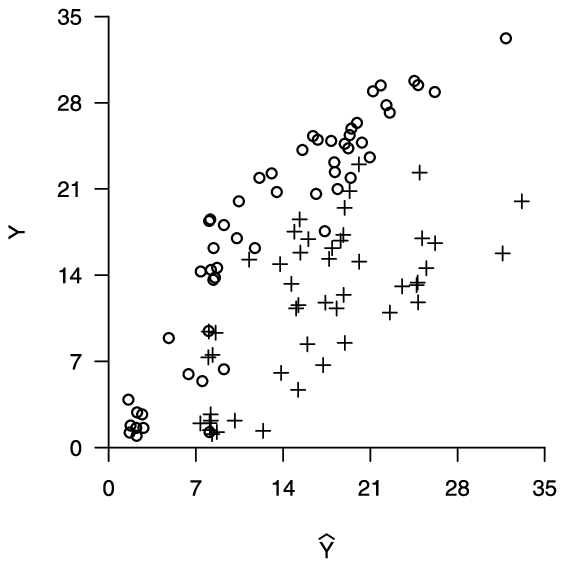}\\
(c) & (d)
\end{tabular}
\caption{Plots of the response versus fitted values for the four
estimators \textup{(a)} $\Rg$, \textup{(b)}~$\Rhatspice$,
\textup{(c)}~$\Rd$ and \textup{(d)} $\Ri$. Circles represent
pork and plus signs represent beef.}\label{figboth}
\end{figure}

Fitting each of the four estimators with $d = 1$, we found that $\rms
(1,3,\allowbreak\Rhat)
=41.96$, $41.65$, $0.55$ and $4.66$ for $\Rhat= \Ri$, $\Rd$, $\Rg$ and
$\Rhatspice$. Plots of~$Y$ versus~$\Yhat$ are shown in
Figure~\ref{figboth}(a)--(d). The relative signal rates for the four
estimators
are reflected by the values of $\Gammabfhat{}^{T}\What\Gammabfhat=
18.4$, $59.2$,
$3346.6$ and $69.3$ for $\What= (p/\tr(\Deltabf))\I_{d}$,
$\diag^{-1}(\Deltabfhat)$, $\Deltabfhat{}^{-1}$ and
$\Deltabfhat{}^{-1}_{\hat\lambda}$.
The relatively large signal for $\Deltabfhat{}^{-1}$ is reflected in the
plots of Figure~\ref{figboth}.

\section{Discussion}\label{secdiscussion}

The class of estimators that we studied is limited in scope relative
to the range of SDR methods presently available for $p = o(n)$ regressions.
However, in the broader context of this article, which does not require
$p =
o(n)$, we introduced the concept of an abundant regression and presented
what may be the first $n,p$ asymptotic analysis of a~class of SDR
methods when the
focus is on estimating a~reduction $\R$ rather than on the underlying
central subspace.
These ideas can in principle be extended for other SDR methods, and we
expect that the same general issues will be encountered.
While each of the methods that we studied can perform usefully in the right
situation, we judge the SPICE and $\Deltabf^{-}$ weight matrices to be the
best overall, although
improvements for nonsparse weight matrices and for regressions with very
highly correlated predictors are still needed.

Our simulation results were all based on normal errors to focus the
presentation and save space. However, we have also conducted a~variety of
parallel simulations using Uniform $(0,1)$, $T_{5}$ and $\chi^{2}_{5}$
errors, each centered and appropriately scaled. The results were essentially
the same as those with normal errors [\citet{CooForRot}].

\subsection{Penalty functions} Alternative penalty functions may be used
for estimating the weight matrix,
particularly in scenarios when the inverse error covariance matrix is not
sparse.
For example, the sparse-seeking penalty function in~(\ref{optomega2}),
$\lambda\sum_{i\neq j} |\Omega_{ij}|$, could
be replaced with the quadratic penalty function,
%
%
\begin{equation} \label{quadpenalty}
\lambda\Biggl( \sum_{i\neq j} \Omega_{ij}^2 + \alpha\sum_{j=1}^{p}
\Omega_{jj}^2 \Biggr),
\end{equation}
where $\alpha\in\{0,1\}$ controls whether or not the diagonal of
$\Omegabf$, the inverse error correlation matrix, is penalized. If
$\alpha=0$, the general SPICE algorithm developed by
\citet{Rotetal08} can efficiently solve~(\ref{optomega2}) with the
penalty defined in (\ref {quadpenalty}). If $\alpha=1$,
\citet{WitTib09} derived an noniterative solution to an equivalent
problem to~(\ref{optomega2}) with the penalty defined in
(\ref{quadpenalty}). Recalling that $\varphi_{j}(\A)$ denotes the $j$th
eigenvalue of a~matrix $\A$, let $\eta_j = \varphi_{j}
(\diag^{-1/2}(\Deltabfhat) \Deltabfhat \diag^{-1/2}(\Deltabfhat))$.
Witten\vspace*{1pt} and Tibshirani showed that the eigenvectors of
$\Omegabfhat_{\lambda}$ are equivalent\vspace*{1pt} to the eigenvectors
of $\diag^{-1/2}(\Deltabfhat)\Deltabfhat\diag^{-1/2}(\Deltabfhat)$ and
that $ \varphi_j ( \Omegabfhat_{\lambda} ) = (4 \lambda)^{-1}\{ (\eta
_{j}^2 + 8\lambda)^{1/2} -\eta_j \}$.

\subsection{Choice of $\f$}

The general rates given in Proposition~\ref{propWhat1} are not very
sensitive to the choice of $\f$ since they hold when $\f$
satisfies the minimal rank condition~(\ref{rank}). Nevertheless, assuming
normality and a~correct $\f$, we obtained the oracle rates of
Proposition~\ref{propWDelta1}, which indicates that
there are advantages to pursuing good choices. The methods sketched
in Section~\ref{choose-f} are often useful in practice, but it is also
possible to develop semiparametric methods to estimate $\bolds{\xi
}$ directly
rather than passing through approximations~$\betabf\f$.
This might be accomplished
iteratively: choose an initial $\f$ and construct the corresponding
estimates $\Gammabf^{1}$, $\bolds{\xi}^{1}=\betabf^{1}\f$ and
$\R^{1}_{\What}$. A new estimate of $\bolds{\xi}$ can be
obtained by
smoothing the coordinates of $\R^{1}_{\What}$ against~$Y$, leading to a~second reduction estimate $\R^{2}_{\What}$. The process can now be
continued until some convergence criterion is met.

\subsection{Variable selection} While we did not incorporate
screening or variable selection into our reduction methodology, the
potential benefits of those procedures are manifested in our results.
Consider,\vspace*{1pt} for instance, a~regression in which $p^{2/3}$ of
the predictors are inactive. Then the oracle rate $\kappa^{-1} =
(p/hn)^{-1/2} = n^{1/2}p^{-2/3}$. However, if we remove $p^{1/3}$ of
the inactive predictors, the oracle rate is increased to $\kappa^{-1} =
n^{1/2}p^{-1/3}$, which should be worthwhile in most applications. Work
along these lines is in progress.

\section*{Acknowledgments}

The authors are grateful to the referees whose
helpful comments led to significant improvements in this article.

\begin{supplement}[id=suppA]
\stitle{Supplement to ``Estimating sufficient reductions of the predictors
in abundant high-dimensional regressions''}
\slink[doi]{10.1214/11-AOS962SUPP} 
\sdatatype{.pdf}
\sfilename{aos962\_supp.pdf}
\sdescription{Owing to space constraints, we have placed the technical
proofs in a~supplemental article [\citet{CooForRot}].
The supplement also contains several preparatory
technical results that may be of interest in their own right and additional
simulations. For instance, we gave in Section~\ref{secsimulations}
simulation results from models with exponentially decreasing error
correlations. In the supplemental article we give parallel results
based on the
same models but with constant error correlations.}
\end{supplement}

%

\printaddresses

\end{document}